\documentclass[12pt,twoside]{article}
\usepackage{amsmath,amsfonts,amssymb,theorem}
 \date{29 Juin 1998}
\title{Invariants de Von Neumann des faisceaux coh\'erents}
\author{Philippe Eyssidieux}

 \newtheorem{prop}{Proposition}[subsection]
 \newtheorem{coro}[prop]{Corollaire}
 \newtheorem{lem}[prop]{Lemme}
 \newtheorem{theo}[prop]{Th\'eoreme}
  
 \newtheorem{exe}[prop]{Exemple}

 \theorembodyfont{\rmfamily}
 \newtheorem{defi}[prop]{D\'efinition}

\newcommand{\qed}{\ifhmode\unskip\nobreak\fi\quad{\ensuremath\square}}

 \newcommand{\hidot}{^{\bullet}}
 
 \newcommand{\wave}{\widetilde}

\begin{document}
 \maketitle

\begin{abstract}

In order to study the group of $L_2$ holomorphic sections of the 
pull-back
to the universal covering space
of an holomorphic vector bundle on a compact complex manifold,
it would be convenient to have a cohomological formalism,
generalizing Atiyah's $L_2$ index theorem.

 In \cite{Eys4}
\cite{Eys5}, such a formalism is proposed in a restricted context.
To each coherent analytic sheaf ${\cal F}$ on a $n$-dimensionnal 
{\em 
smooth projective }
variety
$X^{(n)}$ and each Galois infinite unramified covering
$\pi:\tilde
X \to X$,  whose Galois group is denoted by $\Gamma$,  $L_2$ 
cohomology 
groups denoted by $H^q_2(\tilde X,{\cal F})$ are attached, such that:
\begin{enumerate}
\item The $H^q_2(\tilde X,{\cal
F})$ underly a cohomological functor on the abelian category of 
coherent analytic sheaves on $X$.
\item If ${\cal F}$ is locally free, $H_2^0(\tilde X,{\cal F})$ is 
the group of $L_2$ holomorphic sections of the pull-back
to $\wave X$
of the holomorphic vector bundle underlying ${\cal F}$.
\item  $H^q_2(\tilde X,{\cal
F})$ belongs to a category of  $\Gamma$-modules
on which a dimension function  $\dim_{\Gamma}$ with real values is 
defined.
\item Atiyah's $L_2$ index theorem holds \cite{Ati}:
$$
\sum_{q=0}^{n}(-1)^q \dim_{\Gamma}H^q_2(\tilde X,{\cal
F})=\sum_{q=0}^{n}(-1)^q \dim H^q_2( X,{\cal
F})
$$
\end{enumerate}

The present work constructs such a formalism in the natural 
context of complex analytic spaces.  Here is a  sketch of the main 
ideas of this construction, which is  more or less a Cartan-Serre
version of \cite{Ati}. 

A major ingredient  will be  the construction 
\cite{Farb1}  of an abelian category 
$E_f(\Gamma)$ containing every closed $\Gamma$-submodule of the 
left regular representation. In topology, this device 
enables one to use  familiar homological algebra  
to study 
$L_2$ Betti numbers \cite{Ati}
and Novikov Shubin invariants \cite{NS}. It will play a similar r\^ole here.

 We first construct a $L^p$-cohomology theory ($p\in[1,\infty]$)
for coherent analytic sheaves on a complex space endowed with a 
proper 
action of a group $\Gamma$ such that
 conditions 1-2 are fulfilled. The $L_p$-cohomology on the Galois 
covering 
$\wave  X \to X$ of a coherent analytic sheaf ${\cal F}$ on $X$ is
  the ordinary cohomology of a sheaf on $X$ obtained by an adequate
completion of the tensor product of ${\cal F}$ by the locally 
constant
sheaf on $X$ associated to the left regular representation of the 
discrete group $Gal(\wave X/ X)$ in the space of $L_p$ functions on 
$Gal(\wave X/ X)$.

Then, we introduce an homological algebra device, montelian modules, 
which can be used to calculate the derived category of $E_f(\Gamma)$
and are a good model of the \v Cech complex calculating 
$L_2$-cohomology

Using this we prove that $H^q_2(\tilde X,{\cal F}) \in 
E_f(\Gamma)$, if  $X$ is compact.  This is stronger than
condition 3, since this also yields Novikov-Shubin type invariants. To 
explain the title of the article,
$L_2$ Betti  numbers and
  Novikov-Shubin
invariants of  $H^q_2(\tilde X,{\cal F})$  are the Von Neumann 
invariants 
of the coherent analytic sheaf ${\cal F}$.

We also make the connection with Atiyah's $L_2$-index theorem 
\cite{Ati}
thanks to a Leray-Serre spectral sequence. From this, condition 4
is easily deduced.

\end{abstract}

\section{Introduction}

L'\'etude du groupe des
sections holomorphes $L_2$ du relev\'e au rev\^etement 
universel d'un 
fibr\'e
vectoriel holomorphe sur une vari\'et\'e complexe 
compacte  serait 
grandement
facilit\'ee par un formalisme 
cohomologique\footnote{D\'emarr\'ee par M. 
Gromov \cite{Gro2}, une s\'erie d'investigations r\'ecentes 
utilise cette technique  pour \'etudier les vari\'et\'es 
alg\'ebriques \`a 
groupe fondamental infini et est un des \'el\'ements 
principaux de 
l'id\'eologie de \cite{Kol2}. Parmi les travaux dans cette 
direction, 
citons \cite{Cam}, \cite{Cam2}, \cite{Eys4},\cite{Kol1}, 
 \cite{NR}, \cite{taka}}.

 Dans \cite{Eys4}
\cite{Eys5}, un tel formalisme est propos\'e. A tout 
faisceau
coh\'erent ${\cal F}$ sur une vari\'et\'e complexe {\em projective 
lisse}
$X^{(n)}$ de dimension $n$ et tout rev\^etement galoisien 
non ramifi\'e 
$\pi:\tilde
X \to X$ de groupe de Galois $\Gamma$, on associe des 
groupes de
cohomologie $L_2$ not\'es $H^q_2(\tilde X,{\cal F})$ tels 
que:
\begin{enumerate}
\item Les $H^q_2(\tilde X,{\cal
F})$ s'organisant en un foncteur cohomologique sur la 
cat\'egorie des
faisceaux coh\'erents sur $X$.
\item Si ${\cal F}$ est localement libre, $H_2^0(\tilde 
X,{\cal F})$
s'identifie
au groupe des
sections holomorphes $L_2$ du relev\'e \`a $\tilde X$ du 
fibr\'e
sous jacent \`a ${\cal F}$.
\item Les $H^q_2(\tilde X,{\cal
F})$ appartiennent \`a une cat\'egorie de $\Gamma$-modules
sur laquelle existe une fonction dimension $\dim_{\Gamma}$ 
\`a valeurs 
r\'eelles.
\item On a le th\'eor\`eme d'indice $L_2$ d'Atiyah 
\cite{Ati}:
$$
\sum_{q=0}^{n}(-1)^q \dim_{\Gamma}H^q_2(\tilde X,{\cal
F})=\sum_{q=0}^{n}(-1)^q \dim H^q_2( X,{\cal
F})
$$
\end{enumerate}

Le pr\'esent article a pour but de construire un tel  
formalisme
dans le cadre, naturel pour une telle \'etude, des espaces 
complexes, ce 
qui permet de traiter en 
particulier le cas des vari\'et\'es alg\'ebriques 
singuli\`eres. 

L'\'etude des invariants $L_2$ des vari\'et\'es et complexes
simpliciaux, nombres de Betti $L_2$ (\cite{Ati},
\cite{ChGr},\cite{Dod},\cite{Luc1}) et invariants de 
Novikov-Shubin
(\cite{NS},\cite{NS1} voir aussi 
\cite{GS},\cite{GS2},\cite{LL})
a connu r\'ecemment une simplification remarquable gr\^ace 
\`a M.Farber
\cite{Farb1}, \cite{Farb2} (et ind\'ependamment \`a W. 
L\"uck,
\cite{Luc2})
\footnote{
En g\'eom\'etrie complexe, les nombres de Betti $L_2$ ont 
un comportement
remarquable d\'ecouvert par M. Gromov 
\cite{Gro1},\cite{Gro2} (voir
aussi \cite{ABR}, \cite{Eys2}, \cite{Eys4}, \cite{JZ}).}. 

Le point principal de la contribution  \cite{Farb1} a 
\'et\'e  de construire
une cat\'egorie ab\'elienne $E_f(\Gamma)$ contenant tous 
les sous 
$\Gamma$-modules
ferm\'es de $L_2\Gamma$. Les nombres de Betti $L_2$ et 
invariants de 
Novikov-Shubin
sont des invariants des objets de $E_f(\Gamma)$, ce qui 
permet pour les 
\'etudier
 des arguments
similaires
\`a ceux de  la th\'eorie de l'homologie des complexes
simpliciaux. Dans \cite{Farb2} se trouve  une preuve du
th\'eor\`eme  de Dodziuk \cite{Dod} qui ressemble  \`a la preuve
faisceautique du th\'eor\`eme de De Rham ( cf. par exemple 
\cite{GH}).

Ce gadget d'analyse
fonctionnelle est un ingr\'edient majeur de notre 
construction.

D\'ecrivons l'organisation de cet article.

La section 2 donne les conventions et d\'efinitions de 
base. Nous 
\'elargissons l\'eg\'erement le contexte en autorisant 
que l'action de $\Gamma$ sur $\tilde X$ ait des 
stabilisateurs 
non triviaux et en consid\'erant des faisceaux sur 
lesquels agit une 
extension centrale de $\Gamma$ par $S^1$, de fa\c con \`a 
avoir
une th\'eorie qui inclue le 
\lq Vafa-Witten twisting trick\rq \  de \cite{Gro2}. 

 La section 3 construit
une th\'eorie de la cohomologie $L_p$ des faisceaux 
p\'eriodiques
sur un $\Gamma$-espace complexe, pour $p\in[1,\infty]$, de 
telle
sorte que les conditions 1-2 soient satisfaites.
La cohomologie $L_p$ du rel\'evement \`a un rev\^etement 
galoisien $\wave X \to X$ d'un faisceau coh\'erent ${\cal 
F}$ sur $X$ est la cohomologie ordinaire d'un faisceau sur 
$X$ not\'e $l^p\pi_*{\cal F}$ et obtenu par une compl\'etion ad\'equate du 
produit 
tensoriel de ${\cal F}$
par le faisceau localement constant  sur $X$   associ\'e 
\`a la repr\'esentation r\'eguli\`ere gauche du groupe 
discret $Gal (\wave X/ X)$ sur les fonctions $L_p$ sur
$Gal (\wave X/ X)$. Les difficult\'es  relatives aux 
produit tensoriels compl\'et\'es de faisceaux en espaces 
de Fr\'echet ont \'et\'e contourn\'ees par une 
construction ad hoc. Le prix \`a payer pour l'\'evitement de cette 
difficult\'e est que nous ne savons pas \`a quelle  sous-cat\'egorie 
de la cat\'egorie des  faisceaux analytiques il est possible de prolonger 
le foncteur $l^p\pi_*$.

La section 4 introduit la notion de module hilbertien et 
d\'ecrit
concr\'etement,  
en  paraphrasant \cite{Farb1}, $E_f(\Gamma)$ qui est la 
(plus petite) cat\'egorie ab\'elienne qui  contient la 
cat\'egorie des modules hilbertiens
 puis introduit un outil de
calcul d'alg\'ebre homologique, les modules mont\'eliens, 
qui calculent la cat\'egorie d\'eriv\'ee de $E_f(\Gamma)$ 
et sont un bon mod\'ele pour le
complexe de \v Cech qui calcule la cohomologie $L_2$. 

La section 5
utilise cet outil pour prouver que $H^q_2(\tilde X,{\cal 
F}) \in 
E_f(\Gamma)$, si $X$ est compacte.  Les nombres de Betti 
$L_2$ et les 
invariants de Novikov-Shubin
de $H^q_2(\tilde X,{\cal F})$ sont les invariants de Von 
Neumann du 
faisceau analytique coh\'erent ${\cal F}$. Les invariants 
de Novikov-Shubin d\'etectent le caract\`ere 
non-Hausdorff  de la cohomologie $L_2$. Ces invariants 
sont pour le moment trop mal contr\^ol\'es pour envisager 
des applications en g\'eometrie alg\'ebrique complexe.

La section 6 r\'eduit la condition 4 au  th\'eor\`eme  
d'indice $L_2$
d'Atiyah \cite{Ati} gr\^ace \`a une suite spectrale de 
Leray-Serre.

Les preuves des r\'esultats principaux ont \'et\'e 
obtenues en 
adaptant certaines id\'ees de la preuve du th\'eor\`eme de 
l'image directe 
de Grauert
donn\'ee  par \cite{FoKn} \`a un contexte nettement moins 
compliqu\'e. Les preuves n\'ecessitent souvent un formalisme assez lourd 
et la r\'edaction pr\'esente n'est pas optimale.

Tous ces efforts ne sauraient se justifier que par les
applications. De nouvelles sont en pr\'eparation, mais
de moins nouvelles
\cite{Eys4},\cite{Eys5},
inspir\'ees par le livre de J.Koll\`ar \cite{Kol2},
 utilisent  une version moins g\'en\'erale de la 
pr\'esente th\'eorie.

Signalons aussi un travail ind\'ependant dans une 
direction similaire 
\cite{Cam3}.

Il parait possible de poursuivre le travail ici entrepris 
de 
faisceautisation de la cohomologie 
$L_2$ en direction d'une th\'eorie de Hodge mixte $L_2$, 
peut \^etre d'une 
cohomologie d'intersection $L_2$ pour les espaces 
singuliers. Ceci permet 
probablement d'obtenir des variantes pour des vari\'et\'es singuli\`eres du 
travail de M. Gromov
\cite{Gro2}.

Pour leurs remarques relatives au pr\'esent travail, je 
dois des 
remerciements aux math\'ematiciens suivants: P. Bressler,
F. Campana, J.P. Demailly,
C. Simpson et 
 J. Tapia.

\section{Conventions, Notations}
Le lecteur trouvera dans cette section les notations utilis\'ees librement 
par 
la suite.

\subsection{Divers}

Soit $E$ un espace de Hilbert. ${\cal B}(E)$ d\'esigne 
l'alg\'ebre des 
op\'erateurs born\'es 
de $E$. Soient $E,F$ deux espaces de Hilbert, 
$E\hat\otimes F$ d\'esigne 
leur produit 
tensoriel compl\'et\'e.

Soit ${\cal A}$ une alg\'ebre associative. $Mod_ {{\cal 
A}}$ d\'esigne 
la cat\'egorie des modules sur ${\cal A}$. Si $X$ est un 
espace topologique
$Mod_{{\cal A}}(X)$ d\'esigne la
 cat\'egorie des faisceaux de ${\cal A}$-modules.

\subsection{Alg\'ebre homologique}

Soit $C$ une cat\'egorie et $S$ un syst\`eme multiplicatif 
de morphismes.
$C_S$ d\'esignera la cat\'egorie localis\'ee.

Soit $B$ une cat\'egorie additive.
Conform\'ement aux notations de \cite{Ha2}, I.2,
 $K^i(B)$, $i=b,+,-$, d\'esignera
la cat\'egorie dont les objets sont les complexes 
born\'es, resp.
 born\'es inf\'erieurement, resp.
 born\'es sup\'erieurement, d'objets de la cat\'egorie 
additive $B$ et les 
morphismes
 les classes d'homotopie de morphismes de complexes. 
$C^i(B)$, $i=b,+,-$, 
d\'esignera la 
cat\'egorie
dont les objets sont les complexes et les morphismes les 
morphismes de 
complexes.

 $D^i(A)$ $i=b,+,-$ d\'esignera la
cat\'egorie d\'eriv\'ee born\'ee (resp.
 born\'ee inf\'erieurement, resp.
 born\'ee sup\'erieurement)
de la cat\'egorie ab\'elienne $A$, voir \cite{Verd} 
\cite{Ha2} I.4.

Si $A'$ est une sous cat\'egorie localisante de $A$, 
$D^i_{A'}A$ 
est la cat\'egorie d\'eriv\'ee \`a cohomologie dans $A'$.

La notion de $\delta$-foncteur cohomologique utilis\'ee 
ici est la notion 
d\'efinie en III.1, p.205
dans le  livre  \cite{Har}  et pas celle du Lecture Notes 
\cite{Ha2}. 
Jusqu'\`a la section 5 incluse, le formalisme des 
cat\'egories 
d\'eriv\'ees de J.L. Verdier n'est pas utilis\'e de fa\c con essentielle.

\subsection{Espaces complexes et faisceaux analytiques coh\'erents 
p\'eriodiques}

Pour les
fondations de la th\'eorie des faisceaux analytiques
coh\'erents, nous renvoyons le lecteur \`a \cite{GR}.

\paragraph{Espaces complexes p\'eriodiques}

Soit $\Gamma$ un groupe d\'enombrable discret.

Un {\em $\Gamma$-espace complexe}
est
 un triplet
$(\tilde X,\Gamma,\rho)$ o\`u $\tilde X$ d\'esigne un espace
complexe paracompact s\'epar\'e, et $\rho$ une action  proprement 
discontinue
de $\Gamma$ sur $\tilde X$ par biholomorphismes
\footnote{ Le noyau de $\Gamma\to Biholo(\tilde X)$
peut \^etre non trivial mais sera automatiquement fini}.

Un {\em morphisme de $\Gamma$-espaces complexes }
$(\tilde X,\Gamma,\rho) \to (\tilde Y, \Gamma,\rho)$
est une application holomorphe $f:\tilde X \to \tilde Y$
 qui commute \`a l'action de $\Gamma$.

Une {\em $ \Gamma$-vari\'et\'e complexe}  est
un $\Gamma$-espace complexe r\'egulier.

Un $\Gamma$-espace complexe est dit {\em cocompact} si l'action
de $\Gamma$ est cocompacte.

Un $\Gamma$-espace complexe $\tilde X$
sera dit {\em $\Gamma$-Stein} ssi
il porte une fonction strictement
plurisousharmonique lisse
 \footnote{Une fonction
  sur un espace complexe
sera dite lisse
(resp. plurisousharmonique, resp.
strictement plurisousharmonique),
ssi elle s'exprime localement comme
image r\'eciproque par un plongement d'une fonction lisse 
(resp.plurisousharmonique, resp.
strictement plurisousharmonique )
 sur ${\mathbb C}^N$.}
,
 $\Gamma$-invariante et $\Gamma$-exhaustive
\footnote{ C'est \`a dire  telle que la fonction induite sur $\tilde 
X/\Gamma$
 soit exhaustive} et
 $\tilde X = \Gamma X^0$ o\`u
$X^0$ est une r\'eunion
finie de composantes connexes
stabilis\'ee par un sous groupe
 fini de $\Gamma$.

\paragraph{Faisceaux analytiques
coh\'erents p\'eriodiques}

Soit $\bar \Gamma$ une extension
centrale de $\Gamma$ par un groupe
compact ab\'elien $S$ et $\chi$
 un caract\'ere de $S$.
Dans ce qui suit $\bar \gamma, \bar g, ...$ d\'esignera un
 \'el\'ement  de $\bar \Gamma$
dont l'image dans
 $\Gamma=\bar \Gamma / S$ sera not\'ee $\gamma,g,..$.

Un {\em faisceau analytique coh\'erent
$\bar \Gamma$-p\'eriodique}
sur le $\Gamma$-espace complexe $(\tilde X,\Gamma,\rho)$
est un faisceau analytique coh\'erent
${\cal F}$ sur $\tilde X$
muni d'une action compatible de $\bar \Gamma$\footnote{ Par une action 
compatible de $\bar 
\Gamma$
on entend la donn\'ee pour tout
 $\bar  \gamma \in  \bar \Gamma$
d'un isomorphisme $i(\bar \gamma): {\cal F}
\to \gamma^*{\cal F}$ tel que
$i(\bar \gamma_1\bar \gamma_2)=
\gamma_2^* i(\bar\gamma_1) i(\bar \gamma_2)$ et
pour tout ouvert $U$ de $\tilde X$ la representation
naturelle $i:S\to Aut({\cal F}(U) )$
est continue pour la topologie
induite par la structure
d'espace de Fr\'echet de ${\cal F}(U)$.}.

Un {\em  morphisme de faisceaux analytiques coh\'erents
$\bar \Gamma$- p\'eriodiques} est un
morphisme de faisceaux analytiques coh\'erents qui
commute \`a l'action de $\bar \Gamma$.
Le noyau, le conoyau, l'image et la coimage
au sens des faisceaux sur $\tilde X$
d'un
morphisme de
 faisceaux p\'eriodiques
h\'eritent d'une action
 compatible de $\bar \Gamma$.
La cat\'egorie $C_{\bar \Gamma}(\tilde X)$
des faisceaux analytiques coh\'erents
$\bar \Gamma$-p\'eriodiques
sur un $\Gamma$-espace complexe
est une cat\'egorie ab\'elienne.

Soit ${\cal F}$ un faisceau analytique coh\'erent
$\bar \Gamma$-p\'eriodique.
La formule $\phi_{\chi} = \int _S \chi(z^{-1}) i(z) dz$
d\'efinit un endomorphisme de ${\cal F}$ tel que
$\phi_{\chi} i(z) = \chi(z) \phi_{\chi}$. Le lemme
de Nakayama assure que ${\cal F}$ est
 somme directe $\Gamma$-localement finie
des sous faisceaux $Im(\phi_{\chi'})$,
$\chi'$
parcourant le groupe des caract\`eres de $S$.
Ceci ram\`ene l'\'etude de
$C_{\bar \Gamma}(\tilde X)$
\`a l'\'etude de $C_{\bar  \Gamma, \chi} (\tilde X)$
qui est la sous cat\'egorie de $C_{\bar \Gamma}(\tilde X)$
dont les objets sont
les faisceaux p\'eriodiques pour lesquels $\phi_{\chi}$
est l'identit\'e, i.e.: tels que $S$ agisse sur ${\cal F}$
par multiplication par le caract\`ere $\chi$. De tels faisceaux
seront dits $(\bar \Gamma,\chi)$-p\'eriodiques.
$C_{\bar  \Gamma, \chi} (\tilde X)$
est une cat\'egorie ab\'elienne.

On d\'esigne par $V_{\bar \Gamma, \chi}(\tilde X)$,
(resp. $L_{\bar \Gamma,\chi}(\tilde X)$),
 la sous cat\'egorie additive
(resp. la sous cat\'egorie)
de $C_{\bar \Gamma,\chi}(\tilde X)$
form\'ee des faisceaux
localement libres (resp. inversibles).
Ces deux cat\'egories sont
 naturellement \'equivalentes aux cat\'egories des
fibr\'es vectoriels holomorphes
 (resp. des fibr\'es en droites)
$(\bar \Gamma,\chi)$-p\'eriodiques.

La notion d'action compatible fait sens pour des faisceaux plus 
g\'en\'eraux que les faisceaux coh\'erents. On d\'esigne par 
$Mod_{\bar\Gamma,\chi} (\tilde X)$, la cat\'egorie des $O_{\tilde 
X}$-modules munis d'une $(\bar\Gamma,\chi)$-action compatible.

\paragraph{Une construction locale}

\begin{prop}
\label{condition3}
Soit $U$ un espace $\Gamma$-
Stein. Soit $U'$ un ouvert $\Gamma$-invariant  de $U$ tel que
$\Gamma \backslash U'\subset\subset \Gamma \backslash U$.

Soit ${\cal F}$ un faisceau analytique coh\'erent
$(\bar \Gamma,\chi)$-p\'eriodique d\'efini sur $U$. Il existe un faisceau
analytique coh\'erent
$(\bar \Gamma,\chi)$-p\'eriodique  localement
libre ${\cal
V}$ d\'efini sur $U'$
 et un morphisme surjectif, d\'efini sur $U'$, ${\cal V} \to
{\cal F} \to 0$.

 Soit
 un diagramme de faisceaux
analytiques coh\'erents
$(\bar \Gamma,\chi)$-p\'eriodiques
sur $U$ de la forme:
$$
\begin{array}{ccccc}
 & & {\cal F} & & \\
 &   & \downarrow  & & \\
{\cal W}& \to & {\cal G} & \to &0
\end{array}
$$
la ligne du bas \'etant suppos\'ee exacte.
Il existe un faisceau analytique coh\'erent
$(\bar \Gamma,\chi)$-p\'eriodique  localement
libre not\'e ${\cal
V}$, un \'epimorphisme ${\cal V} \to {\cal F} \to 0$ et
 un morphisme  ${\cal V} \to
{\cal W} $ (d\'efinis sur $U'$) tels que le diagramme suivant commute:
$$
\begin{array}{ccccc}
{\cal V} &\to  & {\cal F} &\to  & 0\\
\downarrow  &   & \downarrow  & & \\
{\cal W}& \to & {\cal G} & \to &0
\end{array}
$$
\end{prop}

{\bf{Preuve:}} Rappelons que $U=\Gamma U^0$ o\`u $U^0$
est une composante connexe stabilis\'ee par un sous groupe fini $\Sigma$ 
de $\Gamma$.
Nous supposons d'abord que $\chi=0$ et que $\Sigma$ est le groupe trivial.
 La proposition  est alors une
cons\'equence du th\'eor\`eme B de Cartan.

Si $\chi\not =0$ et $\Sigma$ est le groupe trivial, il est possible de 
construire
un faisceau inversible $(\bar \Gamma ,\chi)$-p\'eriodique, not\'e
$L_{\chi}$.
Donnons une construction.
Sur l'ensemble
$\bar \Gamma \times U^0 \times {\Bbb C}$
il y a une action $L$ de $\bar \Gamma$ et
une action $R$ de $S$ qui commutent qui se d\'efinissent
 par les formules
$L(\bar \gamma).( \bar g,x,w) =(\bar \gamma \bar g,x,w)$
et $R(z)(\bar g,x,w) =(\bar gz^{-1},x,\chi(z)w)$
Par suite il y a sur
$W:= \bar \Gamma \times U^0 \times {\Bbb C} / R(S)$
une action de $\bar \Gamma$ compatible \`a la
 projection naturelle $W \to \Gamma \times U^0 =
\bar  \Gamma \times U^0 \times{0} / R(S)$.
 De plus, $W \to X$
d\'efinit un fibr\'e
en droites holomorphe sur $\Gamma \times U^0$ muni
d'une $({\bar \Gamma,\chi})$-action compatible.
En utilisant le produit tensoriel par $L_{\chi}$ et son inverse,  on
est ramen\'e au cas pr\'ec\'edent.

 Si $\Sigma\not =0$,
notons que nous avons une application \'etale $\Gamma$-invariante
$\pi: V \to U$ avec $V=\Gamma \times U^0$ un $\Gamma$-espace complexe
sans composante connexe fix\'ee. Les consid\'erations pr\'ec\'edentes
s'appliquent
donc \`a $V$. On descend \`a $U$ en utilisant l'\'epimorphisme
 $\pi_*\pi^*{\cal F} \to {\cal F}$, d\'efini  pour tout
faisceau ab\'elien ${\cal F}$, et fonctoriel en ${\cal F}$.

\section{Cohomologie $L_p$ des faisceaux coh\'erents }

Dans ce qui suit $p$ d\'esigne un nombre
r\'eel
$\in[1,\infty[$.
Pour $p=\infty$, les adaptations \`a fournir sont tr\`es ais\'ees.

\subsection{Construction des
Images directes $L^p$ de faisceaux coh\'erents
p\'eriodiques dans le cas  r\'eduit}

Soit$(\tilde X,\Gamma,\rho)$ un $\Gamma$-espace
complexe r\'eduit, $1\to S\to \bar \Gamma \to \Gamma \to 1$ une extension
centrale de $\Gamma$ et $\chi $ un caract\`ere de $S$.
Notons par $\pi$ l'application continue $\tilde X \to X
=\Gamma \backslash \tilde X $. Convenons de noter pour toute
partie $P$ de $X$,
$\tilde P:= \pi^{-1} (P)$.

\subsubsection{Le cas des faisceaux localement libres}

Si $Z$ est un espace complexe r\'eduit, la topologie de la
convergence uniforme sur les compacts induit une structure
d'espace de Fr\'echet sur l'espace $O(Z)$ des fonctions holomorphes sur
$Z$. De plus, on peut  construire une classe de mesures  sur 
$Z$,
 repr\'esent\'ee par une mesure $dv_Z$ dont le support est $Z$,  et telle
 que, pour tout ouvert $U$ de $Z$
les normes $(\|. \|_{L^p(K,dv_Z)})$,  $K$ parcourant les compacts de
$U$, d\'efinissent  la  structure canonique d'espace de Fr\'echet sur 
$O(U)$.
Si $Z$ est une sous vari\'et\'e ferm\'ee
de ${\mathbb C}^N$ de dimension pure $d$, il suffit de poser $dv_Z 
=i^{reg}_*\omega^d$
o\`u $i^{reg}:Z^{reg}\to Z$ est le plongement ouvert de la partie
r\'eguli\`ere et $\omega$ est la restriction \`a $Z^{reg}$
 de la forme K\"ahlerienne d\'efinie par la m\'etrique
euclidienne. Comme de plus $\tilde X$ est un $\Gamma$-espace complexe
r\'eduit, il est loisible de supposer  $dv_{\tilde X}$
 $\Gamma$-invariante, ce que nous ferons d\'esormais.

Soit ${\cal V} \in V_{\bar \Gamma, \chi} (\tilde X)$.
Fixons une m\'etrique hermitienne $h$, $\bar
\Gamma$-invariante, sur le fibr\'e vectoriel holomorphe
$V$ associ\'e \`a ${\cal V}$.

\begin{defi}
Le faisceau $l^p\pi_*{\cal V}$ est
 le faisceau en groupes ab\'eliens sur $X$ d\'efini par:

$$l^p\pi_* {\cal V}(U)=\{ s\in
H^0 (\pi^{-1} (U), V): \ \forall K
\subset\subset U
 \ \int_{\tilde K} |s|_h^p dv_{\tilde X}< \infty\}$$
\end{defi}

Sricto sensu, la formule pr\'ec\'edente d\'efinit
un pr\'efaisceau, mais
les propri\'et\'es de recollement se v\'erifient ais\'ement. Par ailleurs, 
il est \'egalement
clair
que $l^p\pi_* {\cal V}$ est ind\'ependant du choix de la
m\'etrique.

${\cal V} \mapsto l^p\pi_* {\cal V}$ est un foncteur
de $V_{\bar \Gamma, \chi}(\tilde X)$ vers la cat\'egorie des
faisceaux en groupes ab\'eliens sur $X$.

On dispose d'une inclusion fonctorielle canonique
$l^p\pi_* {\cal V}\buildrel{can ({\cal V})}\over 
{\hookrightarrow}\pi_*{\cal V}$.

\subsubsection{Th\'eor\` eme d'exactitude}

 Soit ${\cal V}, {\cal W}$ des objets de
 $V_{\bar \Gamma,\chi} ( \tilde X)$
 et $\phi: {\cal V} \to {\cal W}$
un morphisme d'image
le faisceau analytique coh\'erent
p\'eriodique  ${\cal I}$.

Soit $\tilde U$ un ouvert
$\Gamma$-Stein  de
$(\tilde X,\Gamma,\rho)$.
Soit $P$ est un fonction strictement
plurisousharmonique lisse $\Gamma$-invariante
 et $\Gamma$-exhaustive d\'efinie
sur  $\tilde U$. Ecrivons
 $\tilde U = \Gamma U^0$ o\`u
$U_0$ est une r\'eunion finie de composantes connexes
stabilis\'ee par un sous groupe fini de $\Gamma$.
Pour $\min P < t  $ on d\'efinit
$\tilde U_t = P^{-1} (] -\infty, t[)$
et $U^0_t =\tilde U_t \cap U^0$.

\begin{lem}\label{frechet}
Fixons $t$ tel que
$\min P <t$. Il existe $t'>t$ et une constante
r\'eelle $ C_{t,t'}>0$,
telle que pour toute
section  $s'$ de ${\cal W}$ d\'efinie sur $U^0$
  et
 localement contenue dans 
${\cal I}$,  il existe une section $s$ de  ${\cal V}$
sur $U^0$ telle que $\phi(s) =s'$ et v\'erifiant:
$$\int_{U^0_{t}} \|s\|^p \le C_{t,t'}\int_{U^0_{t'}}
\|s'\|^p$$
\end{lem}
{\bf { Preuve :} }
Consid\'erons l'espace ${\cal W} (U^0)$
des sections  holomorphes de $W$
 sur $U^0$. Pour $t\in {\mathbb R}$, nous
 d\'efinissons une norme $\| .\|_t$ sur
${\cal W} (U^0)$  par la formule:
$$\| \sigma\|^p_t=\int_{U^0_t} |\sigma|^p$$
Cette collection de normes
 induit sur ${\cal W} (U^0)$  une topologie
 d'espace de Fr\'echet, ind\'ependante de $p$.

\begin{lem}\label{closedness}
L'espace ${\cal I}(U^0) \subset {\cal W} (U^0)$
est  un ferm\'e.
\end{lem}
{\bf {Preuve:}} Voir \cite{GR}, p. 46 \footnote{Il
s'agit d'une application
 standard du lemme de Krull, combin\'e avec l'observation
que la convergence $L^p_{loc}$ force la convergence
des coefficients du d\'eveloppement en s\'erie enti\` ere
 en tout point de $\Delta^n$.}. La topologie induite
sur ${\cal I}(U^0)$ est donc de Fr\'echet.

L'application $\phi: {\cal V} (U^0)  \to
{\cal W} (U^0)$
est une application lin\'eaire continue
car $\|\phi(s')\|_t \le
 (\sup_{x\in U^0_t} ||| \phi(x) |||) \|s'\|_{t}$.
Elle est \`a valeurs dans ${\cal I}(U^0)$
 et $\phi: {\cal V}(U^0) \to  {\cal I}(U^0)$
est surjective
par le Th\'eor\`eme B de Cartan, puisque $U^0$ est de Stein.

Une application continue et surjective
entre espaces de Fr\'echet est ouverte, voir \cite{Yos},
p.75.

De ce fait, l'ouvert
$\Omega_t=\{ s \in  {\cal I}(U^0),
 \ \int_{U^0_t}
|s|^p < 1\}=B_{\|.\|_t}(0,1)$ v\'erifie que $\phi
(\Omega_t)$
est un ouvert de ${\cal I}(U^0)$. Par suite $\exists
t'>0,\epsilon >0$ tel que $\{ s \in {\cal I}(\Delta),
 \ \int_{U_{t'}^0}
|s |^p < \epsilon \}$
 soit contenu dans $\phi(\Omega_t)$. On peut 
 toujours supposer $t'>t$.
Poser $C_{t,t'} = 2/ \epsilon$ \'etablit
le lemme \ref{frechet}.

Le lemme \ref{frechet} se
 traduit imm\'ediatement en le:

\begin{coro} \label{hophop}
L'image de
$l^p\pi_*\phi:
l^p\pi_* {\cal V}\to l^p\pi_*{\cal W}$
s'identifie
 sous l'inclusion naturelle $l^p\pi_*{\cal W} \to
\pi_*{\cal W}$
 au sous faisceau  $l^p\pi_*{\cal W}
  \cap \pi_*{\cal
I} $.
\end{coro}
{{\bf Preuve :}} Puisque $\pi$ est un morphisme de Stein,  $\pi_*{\cal I}$ 
est l'image de $\pi_*\phi: \pi_*{\cal V} \to \pi_*{\cal W}$. Par suite, 
$Im(l^p\pi_*\phi)
\subset l^p\pi_*{\cal W} \cap \pi_*{\cal I}$.

Prouvons la r\'eciproque. Un \'el\'ement $\xi_x$ de
 $(\pi_*{\cal I}\cap l^p\pi_*{\cal W})_x$
 est le germe
en $x$ d'une section de ce faisceau
d\'efinie sur un
ouvert de la forme
$\Gamma \backslash\tilde U $ avec $\tilde U$
un ouvert $\Gamma$-Stein, c'est \`a  dire
 par la donn\'ee d'une
collection
 $(f'_{\gamma})_{\gamma \in \Gamma/ \Sigma}$
de section holomorphes
sur $\gamma.U^0$ v\'erifiant pour toute constante r\'eelle
$t>P(x)$, $\sum_{\gamma}
\int_{\gamma U^0_t} |f'_{\gamma}|^p<\infty$
et $f'_{\gamma} \in {\cal I} (\gamma U^0)$.

Choisissons pour $\gamma \in \Gamma/ \Sigma$
 un \'element
$\bar \gamma \in \bar \Gamma$
 qui s'envoie sur $\gamma$
et posons $s'_{\gamma} (y) := \bar \gamma^{-1}f'_{\gamma}(
\gamma y)$.
$s'_{\gamma}$ est une section
 de ${\cal W}$ sur $U^0$ localement
contenue dans ${\cal I}$.
La proposition \ref{frechet} assure que
l'on peut choisir des constantes  $t'>t>P(x)$
et $C_{t,t'}>0$ ind\'ependantes de $\gamma$ et des sections
holomorphes de ${\cal V}$, $s_{\gamma}$,
 telles que
$\phi(s_{\gamma})=s'_{\gamma}$
 et $\int_{U^0_t} |s_{\gamma}|^p
\le C_{t,t'}\int_{U^0_{t'}} |s'_{\gamma}|^p$.

Posons $f_{\gamma}(y)= \bar \gamma s(\gamma^{-1}y)$.
Ceci d\'efinit une section
de ${\cal V}$ sur $\gamma U^0_t$
tel que
 $\phi (f_{\gamma} )=f'_{\gamma}$
 v\'erifiant de
plus
$\sum_{\gamma} \int_{\gamma U^0_{t}}
|f_{\gamma}|^p \le C_{t,t'}
\sum_{\gamma}
\int_{\gamma U^0_t} |f_{\gamma}'|^p$.

Par suite
$f'=\{f'_{\gamma} \}_{\gamma \in \Gamma /
\Sigma} \in l^p\pi_* {\cal V} ( U_t)$.
 Soit $\eta_x$
son germe en $x$. Par construction,
 $l^p\pi_*\phi( \eta_x)=\xi_x$. $\qed$.

\begin{prop} \label{exact}Soit $ {\cal V}
\buildrel{\phi}\over{\to}{\cal W}
\buildrel{\psi}\over{\to}{\cal X}$  une suite
exacte de $C_{\bar \Gamma,\chi}(\tilde X)$ tels que
${\cal V}$, ${\cal W}$, ${\cal X}$ sont localement libres.
Alors $l^p\pi_*{\cal V}
\buildrel{l^p\pi_*\phi}
\over{\to}l^p\pi_* {\cal W}
\buildrel{l^p\pi_*\psi}
\over{\to}l^p\pi_*{\cal X}$
 est \'egalement exacte.
\end{prop}
{\bf Preuve:} En effet $\ker (l^p\pi_*\psi)=
l^p\pi_*{\cal V}
\cap
\pi_*\ker (\psi)$. Comme
$\ker (\psi) =\text{Im} (\phi)$
on peut appliquer le corollaire \ref{hophop} pour
conclure $\ker(l^p\pi_*\psi)=
\text{Im}( l^p\pi_*\phi)$.$\qed$

\subsubsection{Images directes $L^p$
de faisceaux coh\'erents p\'eriodiques}

Soit $\tilde U$ un ouvert $\Gamma$-Stein de
$\tilde X$ et
${\cal F}$ un objet de $C_{\bar \Gamma, \chi} (\tilde X)$.

Soit $R^{\hidot} \buildrel{i(R^{\hidot})}\over{\to} {\cal F}$ une 
$2$-pr\'esentation \footnote{Soit $F$ un objet  d'une cat\'egorie 
ab\'elienne.  
Une {\em pr\'esentation} $R^{\hidot}: R^{-1}\buildrel{f}\over{ \to } 
R^{0}\buildrel{p}\over{ \to} F \to 0$ de $F$ est la donn\'ee d'un 
isomorphisme $p$ de $F$ avec le  conoyau de la fl\'eche $f$. Une {\em 
r\'esolution}
de $F$ est la donn\'ee d'une suite exacte index\'ee par $-\mathbb N$ de la 
forme  $R^{\hidot}: \ldots \to R^{-n-1} \to R^{-n} \to \ldots \to R^{0} 
\to F \to 0$. Soit  $k\in \mathbb N$,  une {\em $k$-pr\'esentation} de $F$
est une suite exacte de la forme $R^{\hidot}: R^{-k} \to R^{-k+1}\to \ldots
\to R^0\to F \to 0$. }
par des faisceaux $(\bar \Gamma,\chi)$-p\'eriodiques localement
libres. On
dispose
de l'inclusion naturelle $l^p \pi_* {R^{\hidot} }
 \buildrel{can(R^{\hidot})}\over{\to} \pi_*R^{\hidot}$
et du morphisme canonique
 $\pi_* R^0\buildrel{\pi_*i(R^{\hidot})}\over{\to}\pi_*{\cal
F}$.

Notons $\Phi (R^{\hidot})$ le sous faisceau de $\pi_*{\cal F}$
d\'efini comme l'image du morphisme $
l^p\pi_*R^0
\buildrel { \pi_*i(R^{\hidot}) can(i(R^{\hidot}))} \over{\to }
\pi_*{\cal F}$.

\begin{lem} Soient $R^{\hidot}_1\buildrel{i_1} \over{\to } {\cal F}$
 et $R^{\hidot}_2\buildrel{i_2} \over{\to } {\cal F}$
deux $2$-pr\'esentations localement libres
de ${\cal F}$.
S'il existe un morphisme de $2$-pr\'esentations $R^{\hidot}_2\to 
R^{\hidot}_1$  
tel que
chaque fl\'eche $R^n_2\to R^n_1$ soit un \'epimorphisme, 
les deux sous faisceaux  $\Phi(R^{\hidot}_1)$ et $\Phi(R^{\hidot}_2)$ de 
$\pi_*{\cal F}$ coincident .
\end{lem}
{\bf {Preuve: }} Le carr\'e commutatif
$$
\begin{array}{ccc}
R^{\hidot}_2 &\buildrel{i_2} \over{\to } &{\cal F}_U \\
\downarrow \phi_2^1& & \downarrow Id \\
R^{\hidot}_1 &\buildrel{i_1} \over{\to } &{\cal F}_U
\end{array}
$$
 induit un carr\'e commutatif:

$$
\begin{array}{ccc}
{ H}^0 (l^p\pi_* R^{\hidot}_2 )
&\buildrel{c(R^{\hidot}_2)} \over{\to } &\pi_*{\cal F}'_U \\
\downarrow \psi_2^1& & \downarrow Id \\
{ H}^0 (l^p\pi_* R^{\hidot}_1) &\buildrel{c(R^{\hidot}_1)}
\over{\to } &\pi_*{\cal F}'_U
\end{array}
$$
$\Phi(R^{\hidot}_1)$ et $\Phi(R^{\hidot}_2)$ sont les images respectives 
de $c(R^{\hidot}_1)$ et $c(R^{\hidot}_2)$.

Le complexe noyau $K^{\hidot}$ de $\phi_2^1$ est un complexe de
faisceaux coh\'erents
localement libres. Il est acyclique en degr\'es $-1,0$ par la suite exacte 
longue associ\'ee \`a la suite exacte courte associ\'ee 
\`a $0\to K^{\hidot}\to R^{\hidot}_2 \to R^{\hidot}_1$.
Par suite $l^p\pi_* K^{\hidot}$
est un complexe acyclique en degr\'es $-1, 0$ par la proposition
\ref{exact}.

 Toujours par la proposition \ref{exact},
on a une suite exacte courte:
 $$0\to l^p\pi_* K^{\hidot} \to l^p\pi_* R^{\hidot}_2
\buildrel{\phi_2^1}\over{\to} l^p\pi_* R^{\hidot}_1 \to 0$$
 La suite exacte longue associ\'ee 
assure
que ${ H}^0 (l^p\pi_* R^{\hidot}_2)
\buildrel{\psi_2^1} \over{\to}{ H}^0 (l^p\pi_* R^{\hidot}_1)$
est un isomorphisme. Par suite $\Phi(R^{\hidot}_1)= \Phi(R^{\hidot}_2)$. 
$\qed$.

Par le premier point de la proposition \ref{condition3}, pour tout ouvert 
$\Gamma$-invariant  $U'$,
contenu dans un ouvert $\Gamma$-Stein $U$ de $X$ avec
$\Gamma\backslash U'\subset\subset\Gamma\backslash U$
 de
$X$,  la restiction de ${\cal F}$ \`a $U'$ admet une $2$-pr\'esentation 
(et m\^eme une r\'esolution)
localement libre p\'eriodique.

De plus, \'etant donn\'ees deux $2$-pr\'esentations
localement libres sur $U'$ un ouvert  $\Gamma$-Stein, pour tout $U''$ 
ouvert $\Gamma$-invariant
de $U'$ tel que $\Gamma\backslash U''\subset\subset\Gamma\backslash
U'$, un raisonnement ais\'e appuy\'e sur le deuxi\`eme point de la 
proposition
\ref{condition3} assure l'existence
d'une troisi\`eme $2$-pr\'esentation localement libre p\'eriodique 
$R^{\hidot}_3$ d\'efinie sur $U''$
avec des morphismes de comparaison $R^{\hidot}_1 \leftarrow R^{\hidot}_3 
\rightarrow R^{\hidot}_2
$, chaque fl\'eche \'etant une surjection.

Il est donc l\'egitime de poser:

\begin{defi}  $l^p\pi_*{\cal F}$
est l'unique  sous faisceau de $\pi_*{\cal F}$
tel que, $R^{\hidot}_U\to {\cal F}|_U$ d\'esignant une $2$-pr\'esentation 
locale
localement libre p\'eriodique de ${\cal F}$ sur un ouvert
$\Gamma$-invariant $U$,
on ait $l^p\pi_*{\cal F}|_U=\Phi(R^{\hidot}_U)$.
\end{defi}

Si 
$\phi:{\cal F} \to {\cal G}$
est un morphisme de faisceaux coh\'erents p\'eriodiques sur $\tilde X$, le 
morphisme
 $\pi_*\phi: \pi_*{\cal F} \to\pi_*{\cal G}$ envoie $l^p\pi_*{\cal F}$
 dans $l^p\pi_*{\cal G}$. Nous d\'efinissons $l^p\pi_*\phi$ comme la 
restriction de $\pi_*\phi$ \`a $l^p\pi_*{\mathcal F}$.

\begin{prop}\label{l2piet}
Le foncteur $l^p\pi_*: C_{\bar \Gamma,\chi} (\tilde X) \to
Mod_{{\mathbb Z}}(X)$,
${\cal F} \mapsto l^2\pi_*{\cal F}$, $\phi\mapsto l^p\pi_*\phi$ est un 
foncteur exact.
\end{prop}
{\bf{Preuve :}} Seul reste \`a v\'erifier que le foncteur
est exact.  La question est locale.

 Soit $R^{\hidot}_{\cal F} \to
{\cal F}$
une r\'esolution localement libre  d\'efinie
sur $\tilde U$ par des
faisceaux p\'eriodiques. Le morphisme naturel ${ H}^0 (l^p\pi_*{\cal
F})_{\tilde U}
 \to l^p\pi_*{\cal F}_U$ est un isomorphisme. En effet, le
noyau
de l'application naturelle $l^p\pi_* R^0_{\cal F} \to
l^p\pi_*{\cal F}$
coincide avec le noyau de $l^p\pi_* R^0_{\cal F} \to \pi_*
{\cal F}$
qui est exactement $l^p\pi_* R^0_{\cal F} \cap \pi_*{\cal
K}$
o\`u ${\cal K}$ est le noyau de $R^0_{\cal F} \to {\cal
F}$.
Le corollaire \ref{hophop} identifie
$l^p\pi_* R^0_{\cal F} \cap \pi_*{\cal K}$
avec l'image de $l^p\pi_*{R^{\hidot}_1}$,
 ce qui prouve l'isomorphisme annonc\'e.

Soit
 ${\cal S}: 0\to {\cal F}'\to {\cal F} \to {\cal F}'' \to 0$
une suite exacte courte dans $C_{\bar \Gamma, \chi}(\tilde X)$. Pour tout
ouvert
$\Gamma$-invariant assez petit $\tilde U$, utilisant \ref{condition3}, on 
peut 
construire
  une suite exacte d\'efinie sur $\tilde U$ de r\'esolutions
localement libres
 $R^{\hidot}_{\cal S} :0\to R^{\hidot}_{\cal F'}\to
 R^{\hidot}_{\cal F} \to R^{\hidot}_{\cal F''} \to 0$ induisant ${\cal S}$ 
en cohomologie.

Gr\^ace \`a \ref{exact}, il vient une suite exacte courte:
$$l^p\pi_* R^{\hidot}_{\cal S}:0\to l^p\pi_*R^{\hidot}_{\cal F'}\to
 l^p\pi_*R^{\hidot}_{\cal F} \to l^p\pi_* R^{\hidot}_{\cal F''} \to 0$$

Toujours gr\^ace \`a \ref{exact}, pour $q\not =0$,
$${ H}^q (l^p\pi_*R^{\hidot}_{\cal F'})=
{H}^q (l^p\pi_*R^{\hidot}_{\cal F})=
{ H}^q (l^p\pi_*R^{\hidot}_{\cal F''})=0$$

La suite exacte longue associ\'ee \`a $l^p\pi_* R^{\hidot}_{\cal
S}$
se r\'eduit \`a la  suite exacte courte
 $$0 \to { H}^0 (l^p\pi_*R^{\hidot}_{\cal F'})\to
{ H}^0 (l^p\pi_*R^{\hidot}_{\cal F'})\to
{H}^0 (l^p\pi_*R^{\hidot}_{\cal F''})\to 0$$

La $0$-suite associ\'ee \`a ${\cal S}$,  
$0\to l^p\pi_*{\cal F'}
\to l^p\pi_*{\cal F}\to l^p\pi_*{\cal F''}\to 0$
est naturellement isomorphe \`a la pr\'ec\'edente et par suite
exacte.$\qed$.

\subsection{Construction des images directes $L^p$ dans le cas g\'en\'eral}

\begin{lem} \label{modif0}
Soit $\tilde Z$ un $\Gamma$ espace complexe, $\tilde X$
un $\Gamma$-espace complexe r\'eduit et $\tilde i:Z\to X$ un
plongement ferm\'e $\Gamma$-\'equivariant. Notons $i$
l'application $i:\Gamma \backslash \tilde Z \to \Gamma \backslash \tilde X$
induite par $\tilde i$. Soit ${\cal F}\in C_{\bar \Gamma,\chi}(\tilde
X)$.
\begin{enumerate}
\item
 Le faisceau $l^p\pi_* \tilde i_*{\cal F}$ s'identifie \`a un faisceau de 
la forme
 $i_*S^p(i)$
o\`u $S^p(i)$ d\'esigne un sous faisceau  de $\pi_*{\cal F}$.
\item
Si $Z$ est r\'eduit $S^p(i)=l^p\pi_*{\cal F}$.
\item
Si $Z$ n'est pas r\'eduit, soit $I$ le radical de $O_Z$ et soit $({\cal 
F}^k)_{k\in {\mathbb N}}$ la
filtration $I$-adique de ${\cal F}$.
$\pi_*{\cal F}$  admet alors la filtration $(\pi_*{\cal F}^k)_{k\in 
{\mathbb N}}$.
 Le sous-faisceau $S^p(i)$ de $\pi_*{\cal F}$
v\'erifie  $Gr_{\pi_* {\cal F}^{\hidot}}^k  S^p(i) = l^p\pi_* Gr_{\cal 
F^{\hidot}}^k{\cal F}
$ \footnote{Dans cet \'enonc\'e, on voit $Gr_{\cal F}^k {\cal F}$ comme un 
faisceau analytique coh\'erent sur $Z^{red}$
et en particulier comme un faisceau sous l'espace topologique sous jacent 
\`a $Z$.}.
\item
$S^p(i)$ est ind\'ependant du plongement $i$ choisi.
\end{enumerate}

\end{lem}
{\bf Preuve:} Le point  1. r\'esulte du fait que
$l^p\pi_* \tilde i _* {\cal F} \subset \pi_*\tilde i _*{\cal 
F}=i_*\pi_*{\cal
F}$.

Pour le point 2., supposons d'abord ${\mathcal F}$ localement libre. 
Notons $
{\mathcal F}_Z$ le faisceau localement libre sur $Z$ obtenu par restriction
de $\mathcal F$ \`a $Z$.
Comme les applications ${\cal F}(U)\to {\cal F}_Z(U\cap Z)$ sont 
continues, $U$
d\'esignant un ouvert de $X$, 
$S^p(i)\subset l^p\pi_* {\cal F}_Z$ . L'autre inclusion $l^p\pi_* {\cal 
F}_Z\subset S^p(i)$
r\'esulte du th\'eor\`eme  de l'image ouverte appliqu\'e aux
morphismes surjectifs d'espaces de Fr\'echet de la forme ${\cal F}_X(U)\to 
{\cal F}_Z(U\cap Z)$, $U$
d\'esignant un ouvert de Stein de $X$. Le cas g\'en\'eral se
ram\`ene au cas localement libre en utilisant des pr\'esentations locales 
localement libres.

Passons au point 3. Le faisceau $\tilde i _* {\cal F}$ admet la filtration 
$I$-adique
$(\tilde i_*{\cal F}^k)_{k\in {\mathbb N}}$.
La suite des gradu\'es de cette  filtration
est $ (Gr_{\pi_*\tilde i\cal F^{\hidot}}^k{\tilde i_*\cal F})_{k\in {\mathbb
N}}$.
Il r\'esulte de \ref{l2piet} que
$Gr_{\pi_* \tilde i{\cal F}^{\hidot}}^k  i_*S^p(i) =
l^p\pi_* Gr_{\pi_*\tilde i\cal F^{\hidot}}^k{\tilde i_*\cal F}$. Ceci en
conjonction avec les points 1. et 2. ach\`eve d'\'etablir  le point 3.

Soient $i,j$ deux plongements dans un espace r\'eduit.
On peut toujours se ramener au cas o\`u ${\cal F}$ est localement libre.
Quitte \`a remplacer $j$ par  le produit $ k=i\times j$, il est
loisible de supposer que $i= p\circ j$ o\`u $p$ est holomorphe.
Cette factorisation implique que $S^p(i) \subset S^p(j)\subset \pi_*{\cal
F}$. S'il existe $N>0$ tel que $I^N=0$, le point 3. implique 
$S^p(i)=S^p(j)$. Le probl\`eme \'etant local, ceci \'etablit le point  
4.$\qed$

Il est donc l\'egitime de poser:

\begin{defi}  \label{modif}$l^p\pi_*{\cal F}$
est l'unique  sous faisceau de $\pi_*{\cal F}$
tel que, pour tout ouvert
$\Gamma$-invariant $ \tilde U$ et tout $\Gamma$-plongement $\tilde i$
 vers un espace complexe r\'eduit $X$,
on ait $l^p\pi_*{\cal F}|_U=S^p(i)$.
\end{defi}

Cette construction est fonctorielle et il r\'esulte de la proposition
\ref{l2piet} que:

\begin{prop}\label{lpexact}
Le foncteur $C_{\bar \Gamma,\chi} (\tilde X) \to
Mod_{{\mathbb Z}}(X)$,
${\cal F} \mapsto l^p\pi_*{\cal F}$ est un foncteur exact.
\end{prop}

\subsection{ Cohomologie $L^p$ d'un faisceau analytique
coh\'erent
p\'eriodique}

Gr\^ace \`a la proposition \ref{lpexact},
il devient l\'egitime de poser la:

\begin{defi}\label{etvoila}
Soit ${\cal F}$ un faisceau analytique coh\'erent
$(\bar \Gamma, \chi)$-p\'eriodique
 sur le $\Gamma$-espace complexe $\tilde X$,
on appelle $q$-i\`eme groupe de cohomologie $L_p$
de ${\cal F}$ le groupe ab\'elien
 $H^q_{(p)} ( \tilde X, {\cal F})
:=H^q(X, l^p\pi_*{\cal F})$.
\end{defi}

Ceci d\'efinit bien s\^ur un foncteur.
Les homomorphismes de Bockstein de la cat\'egorie
des faisceaux en groupes ab\'eliens sur $X$
donnent \'egalement des homomorphismes de Bockstein
$\delta$ et
on a la:

\begin{prop}\label{defonct}
Le syst\`eme de foncteurs $( \{H^q_{(p)}(\tilde X, .)\}_{q\ge
0}, \delta)$
d\'efinit un $\delta$-foncteur cohomologique de
$C_{\bar \Gamma, \chi} (\tilde X)$
 dans la cat\'egorie des groupes ab\'eliens.
\end{prop}

\section{Modules hilbertiens sur une Alg\`ebre de Von
Neumann}
\subsection{Alg\'ebres de Von Neumann finies}

\begin{defi}
On appelle alg\'ebre de Von Neumann (en abr\'eg\'e AVN) finie 
probabilis\'ee un
triplet
$({\cal A},*,\tau)$ o\`u $\cal A$ est une $\Bbb C$-alg\`ebre
unitaire munie d'une involution $*$ antilin\'eaire et d'
une forme lin\'eaire trace $\tau$ v\'erifiant
\begin{itemize}
\item $\tau(a^*)= \overline{\tau(a)}$
\item $\tau(ab)=\tau(ba), \ \tau(Id)=1$
\item $\tau (a^*a) \ge 0, \ \tau (a^*a) = 0
\Longleftrightarrow a=0$
\item $<.,.>: {\cal A} \times {\cal A} \to {\Bbb C}, \
<a,b>=\tau(ab^*)$
est un produit scalaire pr\'ehilbertien sur $\cal A$. Son
compl\'et\'e ${\cal A}_2$ est un espace de Hilbert
s\'eparable
et le morphisme ${\cal A} \to {\cal B} ({\cal A}_2)$
induit
par
la multiplication \`a gauche a une image ferm\'ee pour
la topologie forte des op\'erateurs\footnote{Ou encore, $\cal A$
coincide avec son bicommutant dans cette
repr\'esentation}. 
\item Si $(a_i)$ est une suite major\'ee croissante\footnote{ Un 
\'el\'ement d'un C*-alg\'ebre  est dit
positif ssi il est autoadjoint $(a=a^*)$ et peut
s'\'ecrire sous la forme
$a=bb^*$. C'est \'equivalent au fait que pour toute
*-repr\'esentation de l'alg\`ebre dans un espace de
Hilbert $H$ ,  $\forall x \in H \ \  (ax,x)\ge 0$. }d'\'el\'ements de
$\cal A$ et $a=\sup_i (a_i)$, alors $\tau(a)=\sup_i(\tau(a_i))$.
\end{itemize}
\end{defi}

Le ${\cal A}$-module ${\cal A}_2$ sera
appel\'e le {\em module standard de $({\cal A}, \tau)$}.

Le commutant de $\cal A$ dans ${\cal B }({\cal
A}_2)$
 est d\'ecrit
par l'action \`a droite de l'alg\`ebre oppos\'ee ${\cal A
}^{op}$. $({\cal A}^{op},*,\tau)$ est \'egalement une AVN finie
probabilis\'ee de module standard ${\cal A}_2$.

\subsection{Exemples d'alg\'ebres de Von Neumann finies}

\subsubsection{Alg\'ebres de Von Neumann d'un groupe discret}

 Soit $\Gamma$ un groupe discret.

 On d\'efinit une alg\'ebre involutive \`a  trace $({\Bbb C }
\Gamma , *, \tau)$ en posant $  (\sum a_g g )^*=\sum \bar a_g
g^{-1}, \
 \tau (\sum a_g g )=a_e$.

 Soit
${\Bbb C} \Gamma \buildrel i \over \to {\cal
B}(l^2\Gamma)$
la repr\'esentation r\'eguli\`ere gauche. Soit
$W_l^*\Gamma$ la
sous alg\`ebre de ${\cal B}(l^2\Gamma)$ qui est
l'adh\'erence de
l'image de $i$ pour la topologie forte des op\'erateurs. 
Si $a\in {\Bbb C} \Gamma, \ a_g=(i(a)
\delta
_e,\delta_g)$.
 De ce fait, tout $a \in W_l^*\Gamma $
  se
repr\'esente comme somme  d'une famille
sommable (pour la topologie forte des op\'erateurs) $
a=\sum a_g
i(g)$. On peut
construire une trace $\tau$ sur l'alg\`ebre
 ${\cal A}=W_l^* \Gamma$ par la formule
  $\tau(a)=a_e=(a\delta_e,\delta_e)$.
On a $\tau (aa^*)=\sum |a_g |^2$. Par suite:

\begin{exe}
 $(W_l^*\Gamma,*,\tau)$
est une alg\`ebre de Von Neumann
probabilis\'ee de module standard $l^2 \Gamma$.
\end{exe}

La m\^eme construction en partant  de la repr\'esentation
r\'eguli\`ere
droite fournit $W^*_r \Gamma$ qui est
 le commutant de $W_l^*\Gamma$.

\subsubsection{Construction de $W^*_l (\bar\Gamma,\chi)$}

La variante ici expos\'ee de la construction pr\'ec\'edente est une 
formalisation du
 \lq Vafa-Witten twisting trick\rq \  de \cite{Gro2}.

Soit $\Gamma$ un groupe d\'enombrable discret.
Soit $S$ un groupe ab\'elien compact.
Soit $\chi$ un caract\`ere de $A$.
Soit $1\to S \to \bar \Gamma \to \Gamma\to 1$ une
extension centrale de $\Gamma$ par $S$.

$\bar \Gamma$ est un groupe topologique localement compact
qui porte une mesure de Haar $d\bar \gamma$ quin se trouve \^etre 
biinvariante.

Consid\'erons l'alg\`ebre de convolution $A$ de $\bar
\Gamma$,
i.e.: $A=C^0_o(\bar \Gamma,{\Bbb C})$, munie du produit de
convolution:
$$ f.f'(\bar g) = \int _{\bar \Gamma} f( \bar g \bar
\gamma ^{-1}) f'(\bar \gamma) d\bar \gamma $$
et de l'antiinvolution antilineaire $f^*(\bar g)=
\overline{f(\bar g ^{-1})}$. $A$ agit sur $E=L^2(\bar
\Gamma,{\Bbb C})$ par  produit de convolution \`a gauche.

Soit $E_{\chi}$ le sous espace ferm\'e de $E$ d\'efini
comme l'ensemble des fonctions
$f$
telles que $\forall z \in S, \forall \bar g \in \bar
\Gamma$, on ait $\chi ( z)f(\bar g z)=  f(\bar g)$. Le
projecteur orthogonal  sur $E_{\chi}$ est:
$$(p_{\chi} f )(\bar g)= \int _S \chi (z) f(\bar g z)dz$$

On a la relation de commutation: $ \forall f \in A,
\forall w\in E, \ f. p_{\chi}w = p_{\chi}(
f.w)$. En effet:
\begin{eqnarray*}
(f. p_{\chi } w) (\bar g) &=&
 \int _{\bar \Gamma} \int _S f(\bar g  \bar \gamma ^{-1})
\chi(z) w(\bar \gamma z) d\bar\gamma dz \\
(p_{\chi} (f.w))(\bar g) &=&
\int _S \int_{\bar \Gamma} \chi(z') f( \bar g z' \bar
(\gamma')^{-1} )w(\bar \gamma ') dz'd\bar \gamma' \\
&\buildrel {z=z', \bar \gamma =\bar \gamma '(z')^{-1}}
\over =& \int _{\bar \Gamma}\int _S \chi(z)f(\bar g
zz^{-1} \bar \gamma^{-1} ) w(\bar \gamma z) d\bar\gamma dz
\end{eqnarray*}

D'o\`u un morphisme d'alg\'ebre involutive $A \to p_{\chi}
A p_{\chi} \subset {\cal B} (E_{\chi})$.
On d\'efinit $W^*_l(\bar \Gamma, \chi)$ comme le bicommutant
de $p_{\chi}Ap_{\chi}$ dans ${\cal B}(E_{\chi})$.

 Posons
$\tau_{\chi}(f) = (p_{ \chi} f)(e)$. On a :

\begin{eqnarray*}
\tau_{\chi} (f.g) &=&
\int_{\bar \Gamma}\int_S f(\bar \gamma^{-1})g(\bar \gamma z)
\chi(z) dzd\bar \gamma\\
\tau_{\chi} (g.f) &=&
\int_{\bar \Gamma}\int_S g((\bar \gamma')^{-1})f(\bar
\gamma' z')
\chi(z) dz'd\bar \gamma'\\
&=&
\int_{\bar \Gamma}\int_S g(z\bar \gamma )f(\bar \gamma^{-1})
\chi(z) dzd\bar \gamma
\end{eqnarray*}

De sorte que $\tau_{\chi} (f.g)= \tau_{\chi} (g.f)$
car $S$ est central dans $\bar \Gamma$.

Le caract\`ere contragr\'edient $\chi^c$ de $\chi$,
$\chi^c (z)
=\bar \chi(z)$, d\'efinit un \'el\'ement de $E_{\chi}$
et l'on a: $\tau_{\chi} (f) = ( f. \chi^c,\chi^c)$.
Observons \'egalement que $f.\chi^c= p_{\chi}f$. De l\`a
$\tau_{\chi} (f.f^*) =\int_{\bar \Gamma} |p_{\chi}f(\bar
\gamma)|^2 d\bar \gamma $.

R\'esumons cette construction en la:

\begin{prop}
$(W^*_l(\bar \Gamma, \chi),\tau_{\chi})$ est une AVN
probabilis\'ee de module standard isomorphe \`a $E_{\chi}$.
\end{prop}

\subsection{ Modules hilbertiens  sur une AVN finie
probabilis\'ee}

Soit $({\cal A},*,\tau)$ une AVN finie probabilis\'ee.

 \subsubsection{Modules hilbertiens projectifs }
\begin{defi}
Un $\cal A$-module \`a gauche $\cal W$ est dit  hilbertien
projectif
s'il admet
un produit scalaire hilbertien $(.,.)$ tel que
$(ax,y)=(x,a^*y)$
et se plonge isom\'etriquement dans le  $\cal A$-module
${\cal A}_2\hat \otimes H$, $H$ \'etant un espace de Hilbert
sur lequel $\cal A$ agit trivialement. $\cal W$ est dit
de
type fini
si $H$ peut \^etre pris de dimension finie. Un morphisme
entre
deux modules hilbertiens est une application lin\'eaire
born\'ee
commutant aux actions respectives de $\cal A$.
\end{defi}

\subsubsection{Non-ab\'elianit\'e de la cat\'egorie des modules 
hilbertiens projectifs}

On d\'efinit trois cat\'egories
$P_f({\cal A}) \subset P_{sep}({\cal A}) \subset P({\cal
A})$ comme
suit:
\begin{itemize}
\item Les objets de $P({\cal A})$
sont exactement les ${\cal A}$-
modules hilbertiens.
\item Les morphismes de $P({\cal A})$ sont exactement
 les applications lin\'eaires continues commutant
\`a l'action de ${\cal A}$, ou encore avec les notations
standard
$Hom_{P({\cal A})} (E,F)= L_{{\cal A}}(E,F)$.
\item $P_f({\cal A})$ (resp. $P_{sep}({\cal A})$)
est la sous cat\'egorie pleine  de $P({\cal A})$ form\'ee
des objets de type fini (resp. r\'ealisables
 dans ${\cal A}_2\otimes H$ o\`u $H$ est un espace
de Hilbert s\'eparable).
\end{itemize}

$P({\cal A})$ est une sous cat\'egorie additive
de la cat\'egorie $Mod_{{\mathcal A}}$.
La notion
cat\'egorique de noyau
d'un $P({\cal A})$-morphisme $f:E\to F$
coincide avec la notion usuelle de noyau.
En revanche, la notion cat\'egorique de conoyau (voir \cite{Lan}, ch. IV,
p. 116)
donne $Coker(f) =(E\to  E/ \overline{ \bar f(E)})$.

$P({\cal A})$ n'est  pas en g\'en\'eral
 une cat\'egorie ab\'elienne:
les morphismes stricts
sont les morphismes d'image ferm\'ee.

\subsubsection {Cat\'egorie ab\'elienne des modules hilbertiens}

Cependant, $P({\cal A})$ poss\`ede la
 propri\'et\'e suivante:

\begin{lem}
Soit $V,W$ deux ${\cal A}$-modules hilbertiens  projectifs et $\phi: V
\to W$
un ${\cal A}$-morphisme surjectif (au sens habituel, c'est
\`a dire
$\phi(V)=W$), alors $\phi$ a une section  ${\cal
A}$-lin\'eaire
continue.
\end{lem}
{\bf{Preuve:}} Fixons une m\'etrique hilbertienne sur $W$
compatible \`a ${\cal A}$ (rappelons que cel\`a signifie
$\forall a,x,y: (ax,y)=(x,a^*y)$, i.e.: que l'adjoint de
l'op\'erateur
multiplication \`a gauche par $a$ est la multiplication
\`a gauche
par le conjugu\'e de $a$). L'orthogonal du noyau de $\phi$
est un sous module hilbertien de $V$ isomorphe via $\phi$
\`a $W$
gr\^ace au th\'eor\`eme de l'image ouverte.

Cette propri\'et\'e, que ne partagent pas les cat\'egories usuelles
 d'espaces vectoriels topologiques localement convexes, a pour cons\'equence
  le fait suivant  dont l'observation est d\^ue \`a
M.S. Farber:

\begin{theo}\cite{Farb1} \cite{Freyd} \label{catabel}
Les cat\'egories $E_f({\cal A}) \subset E_{sep}({\cal A})
\subset E ({\cal A})$ d\'efinies comme suit:
\begin{itemize}
\item Les objets de $E({\cal A})$ sont les triplets
$(E_1,E_2,e)$ o\`u $E_1$ et $E_2$ sont deux
 ${\cal A}$-modules hilbertiens et $e$ une
 application lin\'eaire continue
commutant \`a ${\cal A}$.
\item $Hom_{E({\cal A})} ((E_1,F_2,e),(F_1,F_2,f))$
est l'ensemble
 des classes d'\'equivalence de
 paires d'applications lin\'eaires continues
commuttant \` a ${\cal A}$
$(\phi_1:E_1\to F_1, \phi_2:E_2 \to F_2)$ telles que
$\phi_2 e = f\phi_1$ sous la relation d'\'equivalence
$(\phi_1,\phi_2) \sim ( \phi'_1,\phi'_2) \Leftrightarrow
\exists T\in L_{{\cal A}} (E_2,F_1), \ \phi_2'-\phi_2 =
fT$.
\item $E_f({\cal A})$ (resp. $E_{sep}({\cal A})$)
est la sous cat\'egorie
de $E({\cal A})$ form\'ee des objets $(E_1,E_2,e)$ avec
$E_2$ de type fini  (resp. de type d\'enombrable).
\end{itemize}
sont des cat\'egories ab\'eliennes
\end{theo}

On appellera {\em ${\cal A}$-modules hilbertiens} les objets de $ E({\cal 
A})$.
On notera $P_*(\Gamma) :=P_*(W^*_l \Gamma), \ E_*(\Gamma)
:=E_*(W^*_l \Gamma)$ et $P_*(\bar \Gamma,\chi) :=P_*(W^*_l (\bar
\Gamma,\chi))$, $E_*(\bar \Gamma,\chi) :=E_*(W^*_l (\bar
\Gamma,\chi))$.

\subsubsection{Premi\`eres propri\'et\'es de $E({\cal A})$}

\paragraph{Normalisations des objets et morphismes}

 Un objet $E=(E_1,E_2, e)$ de $E({\cal A})$
sera dit
{\em normal} si et seulement si $e$ est injectif.

\begin{lem} \label{normob} \cite{Farb1}
  Tout objet
$E=(E_1,E_2, e)$ de
$E({\cal A})$ est
isomorphe \`a un objet  normal .
\end{lem}

 Ce  lemme r\'esulte de la propri\'et\'e d'excision \cite{Farb1}:
 soit $(E_1,E_2,e)$ un objet de $E({\cal A})$
et $E\subset E_1$ un sous module hilbertien tel que $eE$
est un sous module ferm\'e de $E_2$. Alors
$(E_1,E_2,e) \simeq (E_1/E,E_2/eE, e)$.

Un repr\'esentant $\phi$ d'un  morphisme $\bar \phi$
dans $E({\cal A})$, $\phi= (\phi_1,\phi_2) $,
 entre
objets normaux,  sera dit {\em normal} si et seulement si
$\phi_1$ est
surjectif.

Soient $E=(E_1,E_2,e),F=(F_1,F_2,f)$ deux objets normaux de
$P({\cal A})$. Soit $\phi=(\phi_1,\phi_2) \in Hom_{P({\cal
A})}(E,F)$,
alors
il existe un objet normal de $P({\cal A})$,
not\'e $E'$, et un isomorphisme $\psi_E : E'\to E$ tels
que le
morphisme $\phi'=\phi \psi_E$ ait un repr\'esentant normal \cite{Farb1}.

\paragraph{Noyaux, Conoyaux}

 Soit $\bar \phi: E\to F$
un morphisme entre objets normaux admettant un
 repr\'esentant normal $\phi=(\phi_1,\phi_2)$.
Soit
$E'=(ker(\phi_1),ker(\phi_2),\bar e)$. Le noyau de  $\bar \phi$
est le 
morphisme $k(\phi)\in Hom_{E({\cal A})}(E',E)$ suivant \cite{Farb1}:
$$
k(\phi) =
 \left(
\begin{array}{ccc}
ker(\phi_1)& \buildrel{\bar e}\over{\to} & ker(\phi_2)\\
k( \phi_1)\downarrow& & k(\phi_2)  \downarrow \\
E_1 & \buildrel{e}\over{\to}  & E_2
\end{array}
\right)
$$

Soit $\bar \phi: E\to F$
un morphisme entre objets normaux admettant un
 repr\'esentant normal $\phi=(\phi_1,\phi_2)$.
Posons $C=(F_1 \oplus
E_2,F_2, f+\phi_2)$. Le morphisme $c=(Id_{F_1}\oplus 0,
Id_{F_2})\in Hom_{E({\cal A})}(F,C)$ est
un conoyau de $\bar \phi$ \cite{Farb1}.

\paragraph{Foncteur d'oubli $E({\cal A}) \to Mod_{{\cal A}}$}

\begin{prop}
La sp\'ecification $O( (E_1,E_2,e))=E_2/ e(E_1)$
d\'efinit un foncteur covariant fid\`ele
$O:E({\cal A}) \to Mod_{{\cal A}}$ et  r\'ealise une \'equivalence de 
cat\'egories
 entre  $E({\cal A})$)
et une sous cat\'egorie ab\'elienne de $Mod_{{\cal A}}$.
\end{prop}
{\bf Preuve:} Un morphisme $f$ dans  $E({\cal A})$
induit clairement
un morphisme $O(f_2): E_2/e(E_1) \to F_2/f(F_1)$, car
 tout repr\'esentant $g$
tel que $g_2 =fT$ v\'erifie $O(g_2)=0$.
Les r\'egles de composition
sont respect\'ees et
on a bien  un foncteur $O: E({\cal A}) \to Mod_{{\cal A}}$.

Le point essentiel est que $Hom_{E({\cal A})}(E,F) \to
Hom_{Mod_{{\cal A}}}(O(E),O(F))$ est injectif.
Par \ref{normob}, il est possible de supposer
$E$ et $F$ normaux. Soit $g=(g_1,g_2)$ un morphisme
tel que $O(g)=0$. Il vient $g_2 (E_2) \subset f(E_1)$.
 En particulier, puisque $f$ est injective,
on peut \'ecrire $g_2 = fT$ o\`u $T$ est une application
lin\'eaire (uniquement d\'etermin\'ee)
 commutant \`a ${\cal A}$.
Reste \` a voir que $T$ est continue. Soit $P_n=(x_n,
Tx_n)_n$
une suite de points du graphe de $T$ convergeant vers
$P=(x,y)$. Comme $f$ est continue, $\lim_n fTx_n = fy$.
Mais, $\lim_n fTx_n= \lim_n g_2 x_n = g_2 x=fTx$.
De l'\'equation $fTx= fy$ on tire $y=Tx$.
$T$ est donc une application lin\'eaire
entre espaces de Banach dont le graphe est ferm\'e, c'est
\`a dire
pr\'ecis\'ement une application lin\'eaire continue.

\begin{lem} \label{isomalg}
Soient $A,B$ deux objets de $E({\cal A})$ et
$\phi\in Hom_{E({\cal A})}(A,B)$ un morphisme. Si
$O(\phi)$
est un isomorphisme alors $\phi$ est un isomorphisme.
\end{lem}
{\bf {Preuve :}} Ceci r\'esulte
du fait que tout  noyau (resp. conoyau) dans $E({\cal A})$
de $\phi$
est envoy\'e par le foncteur d'oubli $O$ dans un
noyau (resp. conoyau) de $O(\phi)$.

\paragraph{Projectifs de $E({\cal A})$}

Le foncteur $p:P({\cal A})\to E ({\cal A})$ et
d\'efini au niveau des objets par
$p(E) = (0,E,0)$ r\'ealise $P({\cal A})$ comme une
sous cat\'egorie pleine de la cat\'egorie des objets
 projectifs de $E({\cal A})$, \'equivalente \`a la cat\'egorie des 
projectifs de $P({\cal A})$ \cite{Farb1}.

\paragraph{Objets de torsion}

Un {\em objet de torsion} $T$ dans $E({\cal A})$
est un objet tel que pour tout projectif $P$ de $E({\cal
A})$, $Hom_{E({\cal A})}(T,P) =0$.

Un objet $T=(T_1,T_2,t)$ est de torsion si et seulement si
$tT_1$ est dense dans $T_2$ \cite{Farb1}.

\paragraph{Suite exacte canonique}

 Le foncteur covariant interne de $E({\cal
A})$ d\'efini au niveau des objets par $P( (X_1,X_2,x) )=
(\overline {xX_1}, X_2,i)$ (resp. $T( (X_1,X_2,x) =
(X_1,\overline {xX_1},x)$) est l'identit\'e sur la sous
cat\'egorie pleine des objets projectifs (resp. de
torsion). On a les relations $P^2 =P$, $T^2=T$,
$PT=TP=0$.
De plus, pour tout objet $X$, on dispose de la suite
exacte
canonique:
$$ 0\to T(X) \to X \to P(X) \to 0$$

\subsubsection{Invariants num\'eriques des objets de
$E({\cal A})$}

\paragraph{${\cal A}$-dimension}

Un $\cal A$-module
hilbertien projectif de type fini $W$ r\'ealis\'e
dans ${\cal A}_2 \otimes {\Bbb C} ^ n$ peut se laisser
d\'ecrire
par son projecteur orthogonal $p =(p_{ij})\in M_n( {\cal
A}^{op})$,
et on pose $\dim_{{\cal A}} (W)= \sum_i \tau (p_{ii})$.
Cette d\'efinition
est ind\'ependante de la r\'ealisation. Cette fonction
dimension
se manipule comme dans l'alg\`ebre lin\'eaire classique,
la
principale
diff\'erence est qu'elle peut prendre
des valeurs r\'eelles positives arbitraires:
\begin{itemize}
\item $\dim_{{\cal A}} W=0 \Longleftrightarrow W=0$
\item $\dim_{{\cal A}} {\cal A}_2 =1$
\item $\alpha \in Mor_{P({\cal A})} (V,W) \Longrightarrow
\dim_{{\cal A}} V = \dim_{{\cal A}} Ker (\alpha) +
\dim_{{\cal
A}}
\overline{ Im (\alpha)}$
\item Si $(W_n)$ est une suite d\'ecroissante de $\cal
A$-
modules
hilbertiens de type fini, $\dim_{{\cal A}} \cap_n W_n=
\lim_n \dim_{{\cal A}} W_n$.
\item Si $(W_n)$ est une suite croissante de sous modules
d'un module
hilbertien de type fini $\dim_{{\cal A}} \overline{ \cup_n
W_n}=
\lim_n \dim_{{\cal A}} W_n$.
\end{itemize}

 Nous appelons {\em ${\cal A}$-dimension} d'un objet
de $X\in P_f({\cal A})$, le nombre r\'eel
$\dim_{\cal A} X:=\dim_{\cal A} P(X)$.

Une propri\'et\'e importante est que pour toute suite exacte
$0\to X' \to X \to X'' \to 0$, on a $\dim_{\cal A}
X = \dim_{\cal A} X' + \dim _{\cal A} X''$.

\paragraph{Invariants de Novikov-Shubin}

Soit $NS$ l'ensemble des germes en z\'ero
 de fonctions d\'efinies sur ${\Bbb R}^+$, croissantes,
positives, nulles en $0$. On appelle $NS_d$ le quotient
de
$NS$ par la relation d'\'equivalence dilatationnelle $\sim
_{d}$ d\'efinie par:
$$F\sim _d G \Longleftrightarrow \exists x_o >0  \
\exists
C,c>0, \ \forall x \le x_o, \  cF(cx) \le G(x) \le
CF(Cx)$$

Soit $A=(A_1,A_2,a)$ un objet de $E_f({\cal A})$ . Fixons
deux m\'etriques hilbertiennes $h_1,h_2$ sur $A_1, A_2$.
Alors on peut construire la famille $\{ E(\lambda)
\}_{\lambda \ge 0}$ des projecteurs
spectraux de $a^*a$. On a :
$$a^*a =E(0) + \int_0^{\| a \| ^2} \lambda dE(\lambda)$$

On pose $F^{A,h_1,h_2} (x^2) =\dim_{{\cal A}} E(x) -
\dim_{{\cal A}} E(0)$.

Il est trivial  que, pour tout autre couple de m\'etriques $h'_1,h'_2$,
 $F^{A,h_1,h_2}\sim_{d}F^{A,h_1',h_2'}$ et on peut donc d\'efinir $F^A$
comme la classe de $\sim_{d}$-\'equivalence  de
$F^{A,h_1,h_2}$.

Moins trivial est le fait suivant \cite{Farb1}:
soient $A$ et $A'$ deux  objets isomorphes de $T_f({\cal A})$. Alors,
 $F^A=F^{A'}\in NS_d$ .

L'{\em invariant de Novikov-Shubin} $NS$ d'un objet $A$ de
$E_f({\cal A})$ est $NS(A) = NS( T(A))$.

\subsection{Calcul de la cat\'egorie d\'eriv\'ee de $E_f({\cal A})$
\`a l'aide des ${\cal A}$-modules mont\'eliens}

\subsubsection{Op\'erateurs ${\cal A}$-compacts}

\begin{defi}
 Soient $W,V$ deux ${\cal A}$-modules hilbertiens projectifs.
$\phi\in {\cal L}_{\cal A} (V,W)$
est dit ${\cal A}$-compact si et seulement si
il existe $(\phi_n)_{n\in {\Bbb N}}$ une suite
d'op\'erateurs convergeant vers $\phi$ au sens de la norme
triple tels que, pour tout $n$,
l'adh\'erence de l'image de $\phi_n$ est
un sous module hilbertien projectif de type fini de $W$.
\end{defi}

Si l'un des deux modules $V$ ou  $W$ est de type fini,
tout \'el\'ement
${\cal L}_{\cal A}(V,W)$ est ${\cal A}$-compact.

Si $\phi \in {\cal L}_{\cal A}(V,W)$ est ${\cal A}$-compact
et surjectif, $W$ est de type fini.

\begin{lem} \label{lemco} Soit $\phi \in {\cal L}_{\cal A}
(V,W)$.
Les
assertions suivantes sont \'equivalentes:
\begin{enumerate}
\item $\phi$ est ${\cal A}$-compact.
\item $\{ E(\lambda)\}_{\lambda \ge 0}$
d\'esignant la famille des projecteurs spectraux de
$\phi^*\phi$, le projecteur $P_{\epsilon}
=\int_{\epsilon}^{\infty} dE(\lambda)$ a pour image un
module hilbertien projectif de type fini pour tout $\epsilon >0$.
\item Pour tout $\epsilon >0$ il existe un sous module
hilbertien projectif $V_{\epsilon}$ de V  de type fini
de projecteur orthogonal $p_{V_{\epsilon}}$, tel que
$\lim_{\epsilon \to 0} \phi p_{V_{\epsilon}}=\phi$.
\item
Pour tout $\epsilon >0$ il existe un sous module
hilbertien projectif $W_{\epsilon}$ de W  de type fini
de projecteur orthogonal $p_{W_{\epsilon}}$, tel que
$\lim_{\epsilon \to 0} p_{W_{\epsilon}}\phi=\phi$.
\end{enumerate}
\end{lem}

\begin{coro} \label{coroco} Le compos\'e (\`a  droite ou \`a  gauche) d'un 
morphisme
${\cal A}$-compact avec une application ${\mathcal A}$-lin\'eaire
continue est ${\cal A}$-compact. En particulier,
la restriction d'un ${\cal A}$-morphisme
${\cal A}$-compact \`a un sous module hilbertien projectif ferm\'e
de sa source est encore ${\cal A}$-compact.
\end{coro}

\subsubsection{ ${\cal A}$-modules mont\'eliens}
\begin{defi}
Soit $t_0\in {\mathbb R}^*_+$.
Un  ${\cal A}$-module $t_0$-pr\'emont\'elien
est une donn\'ee de la forme
 $W=(\{ W(t)
\}_{t_0\ge t>0},\{(\rho^W)_{t}^{t'}\}_{0<t'\le t \le t_0})$
o\`u:
\begin{itemize}
\item $\{W(t)\}_{t_0\ge t>0}$ est une famille de  ${\cal A}$-module
hilbertien projectif index\'ee pour l'intervalle r\'eel
  $]0,t_0]$.
\item Pour tout couple de
 r\'eels tels que $0<t'\le t\le t_0$,
 $(\rho^W)_{t}^{t'}$ est une
application ${\cal A}$-lin\'eaire continue
de $ W(t)$ vers $W(t')$ .
\item $(\rho^W)_{t}^{t} =Id_{W(t)}$.
\item Si $t\ge t'\ge t'' >0$, on a $(\rho^W)_{t}^{t''}=
(\rho^W)^{t''}_{t'} (\rho^W)_{t}^{t'}$
\item $(\rho^W)_{t}^{t'}$ est ${\cal A}$-compacte d\`es que $t'<t$.
\end{itemize}

\end{defi}

\begin{defi} Soient $V,W$ deux ${\cal A}$-modules
$t_0$-pr\'emont\'eliens. Une  famille de la forme
$\{ \phi(t) \}_{t_0\ge t>0}$
de  ${\cal A}$-applications lin\'eaires
continues
$\phi(t):V(t)\to W(t)$ est dite un morphisme de modules
$t_0$-pr\'emont\'eliens si et seulement si
$\phi(t') (\rho^V)_{t}^{t'} = (\rho^W)_{t}^{t'} \phi(t)$.
\end{defi}

\begin{defi} \label{modulana}
Un module pr\'emont\'elien est un module $t_0$-pr\'emont\'elien pour
un certain $t_0>0$ non sp\'ecifi\'e.

Deux ${\cal A}$-modules
pr\'emont\'elien $W$
$W'$ sont dits $G$-\'equivalents s'il existe $t_1>0$ tel que
si $t\le t_1$ $W(t)=W'(t)$ et si $t'\le t \le t_1$,
$(\rho^W)_t^{t'}=(\rho^{W'})_t^{t'}$.

Deux morphismes
$\phi,\phi'$ de ${\cal A}$-modules pr\'emont\'eliens
sont dits $G$-\'equivalents
si leurs sources (et leurs buts) sont $G$-\'equivalents et que
pour $t>0$ assez petit
$\phi(t)=\phi'(t)$.
\end{defi}

\begin{defi}

Un ${\cal A}$-module mont\'elien est
une classe de $G$-\'equivalence de
${\cal A}$-modules pr\'emont\'eliens.
Un morphisme de ${\cal A}$-modules mont\'eliens
est une classe de $G$-\'equivalence de morphismes
de ${\cal A}$-modules pr\'emont\'eliens.
\end{defi}

La cat\'egorie des ${\cal A}$-modules mont\'eliens
est  une cat\'egorie additive o\`u les noyaux existent
(par \ref{coroco}).

Soit $W$ un ${\cal A}$-module pr\'emont\'elien.
Pour ne pas alourdir les notations, on sous-entendra d\'esormais
l'exposant
$W$ et on \'ecrira
$(\rho^W)^{t'}_t =\rho_t^{t'}$.

Pour $t>0$ assez petit $( \{W(t') \}_{t'<t}, \rho_{t''}^{t'})$
forme un syst\`eme projectif
 index\'e par l'ensemble ordonn\'e $(] 0, t[, >)$.

\begin{defi}
$W_*(t)=\underleftarrow{\lim} _{t'<t} W(t')$ 
\end{defi}

$W_*(t)$ est un ${\cal A}$-module de Fr\'echet, dont la
topologie peut \^etre  d\'efinie par une collection d\'enombrable
de pseudonormes
pr\'ehilbertiennes compatibles.

 $\{ \rho_{t}^{t'}\}$ d\'efinissent des applications ${\cal
A}$-lin\'eaires
continues
dans les situations suivantes:

\begin{itemize}
\item $t\ge t'$, $\rho_{t}^{t'}: W(t) \to W(t')$.
\item $t\ge t'$, $\rho_{*t}^{*t'}: W_*(t) \to W_*(t')$.
\item $t\ge t'$, $\rho_{t}^{*t'}: W(t) \to W_*(t')$.
\item $t> t'$, $\rho_{*t}^{t'}: W_*(t) \to W(t')$.
\end{itemize}

On a, de plus, quand les applications en question sont
d\'efinies, $\rho_{a}^{b}\rho_{c}^a= \rho_{c}^b$, $a,b,c$ \'etant
des symboles de la forme $t$ ou $*t$, $t\in{\Bbb R}$.

\subsubsection{ Qi-morphismes}

Soit $(F^{\hidot},d)$ un complexe  de ${\cal A}$-modules
mont\'eliens. Le module des cycles $Z^q(F^{\hidot})$ est un module
mont\'elien, mais pas n\'ecessairement celui des bords.

On dispose d'une  famille de complexes de ${\cal A}$-modules
$(F_*(t),d)_{t>0}$ d\'efinie pour $t>0$ assez petit et l'on pose:

\begin{eqnarray*}
Z^i_*(t,F^{\hidot}) &=& \{ z \in F^{i}_* (t); dz=0\} \\
B^i_*(t,F^{\hidot}) &=& d F^{i-1}_* (t) \\
H_*^i(t, F^{\hidot}) &=& Z_*^i(t,F^{\hidot}) / B_*^i (t,F^{\hidot})
\end{eqnarray*}

$H_*^i(t,F^{\hidot})$ a une structure topologique de quotient non
s\'epar\'e d'espace de Fr\'echet
dont nous ne nous soucierons pas en nous contentant de
le regarder que comme un ${\cal A}$-module au sens alg\'ebrique.

L'identit\'e du complexe $F^{\hidot}$ induit pour
$t'\le t$ des
applications ${\cal A}$-lin\'eaires $ \rho_{*t}^{*t'}: H_*^i(t,F^{\hidot}) 
\to
H_*^i(t',F^{\hidot})$.

Tout morphisme de complexes de ${\cal A}$-modules
mont\'eliens $\sigma ^{\hidot}: F^{\hidot} \to G^{\hidot}$ induit donc des
applications ${\cal A}$-lin\'eaires $ (\sigma^i)_t^{t'}:
H_*^i(t,F^{\hidot}) \to H_*^i(t',G^{\hidot})$, d\'efinies par la formule
$(\sigma^i)_{t}^{t'}= \rho_{*t}^{*t'} H^i(\sigma)_{*t}^{*t}$.

\begin{defi} $\sigma^{\hidot}: F^{\hidot} \to G^{\hidot}$ sera appel\'e un
qi-morphisme d'ordre $q\in {\Bbb Z}$ si et seulement
si pour
tous $t'\le t$ assez petits
\begin{itemize}
\item  $i>q, t'\le t \Rightarrow (\sigma ^i)_t^{t'}
:H_*^i(t,F^{\hidot}) \to H_*^i (t',G^{\hidot})$ est un isomorphisme
alg\'ebrique.
\item  $t'\le t \Rightarrow ( \sigma ^q)_t^{t'}
:H_*^q(t,F^{\hidot}) \to H_*^q (t',G^{\hidot})$ est une surjection.
\end{itemize}
$\sigma^{\hidot}$ est dit  un qi-morphisme si et seulement
si c'est un qi-morphisme de tout ordre.
\end{defi}

\subsubsection{Complexes qhtf}

Un complexe $L^{\hidot}$ de ${\cal A}$-modules hilbertiens projectifs de 
type fini
peut \^etre vu comme un complexe de
${\cal A}$-modules mont\'eliens en posant $W^{\hidot}(t)=L^{\hidot}$
et $\rho^t_{t'}=Id$. En particulier l'identit\'e d'un
tel complexe est un qi-morphisme, puisque
$H_*^i (t,L^{\hidot})=O(
H^i(L^{\hidot}))$ (cf. les notations du  lemme \ref{isomalg}).

\begin{defi}
Un complexe $F^{\hidot}$ de modules
 mont\'eliens sera dit qhtf
si et seulement si
 $Id:F^{\hidot}\to F^{\hidot}$ est un qi-morphisme.
\end{defi}

Soit $\sigma^{\hidot}:F^{\hidot} \to G^{\hidot}$ un morphisme de complexes
de ${\cal A}$-modules mont\'eliens. Le c\^one de $\sigma^{\hidot}$,
$C^{\hidot}(\sigma)$
est le complexe  dont les espaces de cochaines
sont les $C^i(\sigma) = G^i \oplus F^{i+1}$ et les
diff\'erentielles sont $\partial (g^i,f^{i+1}) = ( dg^i +
\sigma f^{i+1} ,-df^{i+1})$. On a la suite exacte courte
$0 \to G^{\hidot} \to C^{\hidot}(\sigma) \to F^{\hidot}[1]\to 0$ qui 
induit des
suites exactes courtes $0 \to G_*^{\hidot}(t) \to
C_*^{\hidot}(\sigma)(t) \to F_*^{\hidot}[1](t)\to 0$.

\begin{lem}\label{qiq} Soit $\phi:P^{\hidot} \to F^{\hidot}$ un
morphisme de
complexes
de ${\cal A}$-modules mont\'eliens, le complexe $P^{\hidot}$ \'etant un
complexe de ${\cal A}$-modules hilbertiens projectifs de type fini. $\phi$ 
est
un qi-morphisme d'ordre $q$ si et seulement si
$\forall i\ge q$,
$\forall t>0$ assez petit, $H_*^i(t,C^{\hidot}(\sigma))=0$.
\end{lem}
{{\bf Preuve :}} L'implication directe est bien connue, prouvons la 
r\'eciproque. 
 Supposons donc que $\forall t>0$ assez petit, 
$H_*^i(t,C^{\hidot}(\sigma))=0$. Soit $0<t_1\le t_0$ assez petits.
Utilisant la suite exacte longue
du c\^one, on voit que $H_*^i(t_l,P^{\hidot}) \to H_*^i(t_l, L^{\hidot})$
est surjectif pour $i=q$ et un isomorphisme pour $i>q$ et $l=0,1$.
Mais, $ \ H_*^i(t_0,P^{\hidot} ) \to H_*^i(t_1,P^{\hidot})$
est, par d\'efinition, l'identit\'e de $O(H^i(P^{\hidot}))$; donc 
$H^i_*(t_0,P^{\hidot}) \to H^i_*(t_1,L^{\hidot})$ 
est surjectif pour $i=q$ et un isomorphisme pour $i>q$.$\qed$.

Plus g\'en\'eralement, ce lemme vaut quand $P^{\hidot}$ est qhtf.

\begin{lem} \label{five}
Soit $\sigma^{\hidot}: F^{\hidot} \to G^{\hidot}$ un morphisme de complexes
de ${\cal A}$-modules mont\'eliens. Si $F^{\hidot}$ et $G^{\hidot}$ sont 
qhtf
le c\^one $C^{\hidot}(\sigma)$ est \'egalement qhtf.
\end{lem}
{\bf{Preuve :}} Il suffit d'utiliser la suite exacte longue de
cohomologie du c\^one et d'appliquer le lemme de cinq.$\qed$.

\begin{prop}
La sous cat\'egorie pleine  $K^i(MontMod_{\cal A})^{qhtf}$, $i\in\{b,-\}$,
 de la cat\'egorie triangul\'ee $K^i(MontMod_{\cal A})$
 dont les
objets sont les complexes qhtf est une sous cat\'egorie
triangul\'ee.

Les qi-morphismes forment un syst\`eme multiplicatif $Qi$ de la
cat\'egorie triangul\'ee $K^b(MontMod_{\cal A})^{qhtf}$ (resp. de 
$K^-(MontMod_{\cal
A})^{qhtf}$).
\end{prop}
{\bf Preuve:} Le premier point d\'ecoule de \ref{five}.

 Le deuxi\`eme point d\'ecoule 
 de \cite{Ha2}, I.4, proposition 4.2; on prend comme foncteur 
cohomologique sur la cat\'egorie triangul\'ee des complexes
qhtf vers la cat\'egorie des groupes ab\'eliens le foncteur $H:C^{\hidot} 
\mapsto \underrightarrow{\lim}_{t\to 0} H^0(t,C^{\hidot})$
\footnote{Que les qi-morphismes soient exactement les morphismes qui 
deviennent des isomorphismes sous $H$ r\'esulte du fait que le syst\`eme 
dont on prend la limite directe est un syst\`eme dirig\'e 
d'isomorphismes, puisqu'on a pris le soin de se limiter aux complexes 
qhtf.}.

\subsubsection{Th\'eor\`eme de finitude}

\begin{prop} \label{tfini} Soit $F^.$ un complexe
fini (resp. born\'e sup\'erieurement) qhtf de
${\cal A}$-modules
mont\'eliens , alors il existe un complexe fini
(resp. born\'e sup\'erieurement) de modules
hilbertiens projectifs de type fini $P^{\hidot}$ et un qi-morphisme
$P^{\hidot} \to F^{\hidot}$.

De plus, pour tout complexe
fini
(resp. born\'e sup\'erieurement) $Q^{\hidot}$ de modules hilbertiens
projectifs
de type fini
et tout morphisme $\phi ^{\hidot}: Q^{\hidot} \to F^{\hidot}$, il existe
un complexe
fini
(resp. born\'e sup\'erieurement)
de modules hilbertiens projectifs de type fini $P^{\hidot}$,
un qi-morphisme $\sigma^{\hidot}: P^{\hidot} \to F^{\hidot}$ et
un morphisme de complexes $\psi^{\hidot}: Q^{\hidot}\to P^{\hidot}$
tels que $\phi^{\hidot} = \sigma^{\hidot} \phi^{\hidot}$.
\end{prop}

{\bf{Preuve :}} Prouvons d'abord la premi\`ere assertion.
 
Le complexe trivial s'envoie sur $F^{\hidot}$ et donne un qi-morphisme 
d'ordre $Q$ pour $Q$ assez grand. Il suffit
de donc  montrer qu'\`a tout qi-morphisme d'ordre $q$,
$\sigma^{\hidot}:L^{\hidot} \to F^{\hidot}$, $L^{\hidot}$ complexe de 
modules
hilbertiens de type fini, on sait associer un nouveau
complexe de modules hilbertiens de type fini $\hat L^{\hidot}$ et
un qi-morphisme d'ordre $q-1$,
$\hat \sigma^{\hidot}: \hat L^{\hidot} \to F^{\hidot}$ .
\footnote{ Dans le cas d'un complexe fini,
si $q$ est assez petit
assez petit ($F^i=0, \ i\le q$), et $\sigma^{\hidot}:L^{\hidot}\to 
F^{\hidot}$
est un qi-morphisme d'ordre $q$,
 le complexe $L_{\infty }^{\hidot}$ avec $L^q_{\infty}=L^q/Z^q(L^{\hidot})$
$L_{\infty}^i=0, i<q$ et $L^i_{\infty}=L^i, i>q$
et le morphisme $\sigma_{\infty}^{\hidot}$ induit par $\sigma^{\hidot}$
d\'efinissent un qi-morphisme $L_{\infty}^{\hidot} \to F^{\hidot}$}

Soit $(L')^{\hidot}$ un complexe de modules hilbertiens projectifs de type
fini.
Soit $\sigma':(L')^{\hidot}\to F^{\hidot}$ un qi-morphisme d'ordre
$q$.
Le complexe $L=T_{q-1}L'$ obtenu en tronquant en $q-1$ le complexe
$L'$ (ie
$(T_{q-1} L')^i=L^i, i\ge q$ et $L^i= 0 $ si $i<q$) donne \'egalement lieu
\`a
un qi-morphisme d'ordre $q$, $\sigma: L^{\hidot} \to F^{\hidot}$.

\begin{lem} \label{id1} Soient $0<t''<t'<t$ assez petits.
Alors:
$$\rho_{t'}^{t''}Z^{q-1}(C^{\hidot}(\sigma))(t') \subset
\rho_{t}^{t''}Z^{q-1}(C^{\hidot}(\sigma))(t) +
d
F^{q-2}(t'')$$
\end{lem}
{\bf {Preuve:}} $\sigma$ v\'erifie $C^{q-2}(\sigma) =
F^{q-2}, i>1$.

Par le lemme \ref{five}, $C^{\hidot}(\sigma)$
est qhtf. Fixons donc $t_0>0$ assez petit
de sorte que $t_1 \le t_0\Rightarrow H_*^q(t_0,
C^{\hidot}(\sigma))\to
H_*^q(t_1,C^{\hidot}(\sigma))$ est un isomorphisme. On peut supposer
$0<t''<t'<t<t_0$ de sorte que
$Z^{q-1}_*(t',C^{\hidot}(\sigma))\subset
\rho_{*t_0}^{*t'}Z^{q-1}_*(t_0,C^{\hidot}(\sigma))+ \partial C^{q-2}_*(t')$.
Puis:
$$\rho_{*t'}^{t''}Z^{q-1}_*(t',C^{\hidot}(\sigma))\subset
\rho_{*t_0}^{t''}Z^{q-1}_*(t_0,C^{\hidot}(\sigma))+
\rho_{*t'}^{t''}\partial C^{q-2}_*(t')$$

L'inclusion annonc\'ee r\'esulte des inclusions suivantes:
\begin{itemize}
\item
$\rho_{t'}^{t''}Z^{q-1}(C^{\hidot}(\sigma))(t')=
\rho_{*t'}^{t''} \rho_{t'}^{*t'} Z^{q-1}(C^{\hidot}(\sigma))(t')\subset
\rho_{*t'}^{t''} Z^{q-1}_*(t',C^{\hidot}(\sigma))$.
\item
$\rho_{*t_0}^{t''}Z^{q-1}_*(t_0,C^{\hidot}(\sigma))
=\rho_t^{t''}\rho_{*t_0}^{t}Z^{q-1}_*(t_0,C^{\hidot}(\sigma))
\subset \rho_{t}^{t''}Z^{q-1}(C^{\hidot}(\sigma))(t)$
\item
$\rho_{*t'}^{t''} C^{q-2}_*(t')\subset F^{q-2}(t'')$.
\end{itemize}

\begin{lem}
Soient $0<t''<t'<t$ assez petits. Il existe deux morphismes
de modules hilbertiens projectifs $S_{t'}^t: Z^{q-1}(C^{\hidot}(\sigma))(t')
\to
Z^{q-1}(C^{\hidot}(\sigma))(t)$ et
$F_{t'}^{t''}:Z^{q-1}(C^{\hidot}(\sigma))(t')
\to F^{q-2}(t'')$ tels que:
$$\rho_{t'}^{t''}|_{Z^{q-1}(C^{\hidot}(\sigma))(t')} =
\rho_{t'}^{t''}\rho_{t}^{t'}S_{t'}^t +dF_{t'}^{t''}$$

\end{lem}

{\bf{Preuve :}}

On pose $Z_{t',t''}= Z^{q-1}(C^{\hidot}(\sigma))(t') /
\ker(\rho_{t'}^{t''})$.

On note $r_{t'}^{t''}: Z_{t',t''}\to
Z^{q-1}(C^{\hidot}(\sigma))(t'')$ l'application induite par
$\rho_{t'}^{t''}$.

On note $p: Z^{q-1}(C^{\hidot}(\sigma))(t')\to Z_{t',t''}$
l'application de passage au quotient.

 $r_{t'}^{t''}$
est injective et  $r_{t'}^{t''} p =\rho_{t'}^{t''}$.

L'image $I_{t',t''}$ de $\rho_{t'}^{t''}$ est la m\^eme que celle
de $r_{t'}^{t''}$. $r_{t'}^{t''}:Z_{t',t''} \to I_{t',t''}$
poss\'ede un inverse $(r_{t'}^{t''})^{-1}: I_{t',t''} \to
Z_{t',t''}$ qui est   une application ${\cal A}$-lin\'eaire
bien d\'efinie, en g\'en\'eral non continue. Posons:
$$C'= \{ (\zeta, f) \in Z^{q-1}(C^{\hidot}(\sigma)) (t)\oplus
F^{q-2}(t''),
\ \text{tels que} \ \ \rho_{t}^{t''}\zeta + d f \in
\rho_{t'}^{t''} Z^{q-1}(t',
C^{\hidot}(\sigma)) \}$$

Sur ce ${\cal A}$-module, la norme
pr\'ehilbertienne
$$\|(\zeta,f)  \|^2_{C'}=
\| \zeta \|^2_t + \| f \| ^2_{t''} + \|( r_{t'}^{t''})^{-1}
(\rho_t^{t''} \zeta +
df)\|_{t'}^2$$
 est compl\`ete et compatible \`a ${\cal A}$.

L'inclusion naturelle isom\'etrique  $C' \to
Z^{q-1}(C^.(\sigma)) (t)\oplus F^{q-2}(t'') \oplus Z_{t',t''}$
montre que c'est un ${\cal A}$-module hilbertien.

 Gr\^ace au lemme
\ref{id1}, il vient:
$$r_{t'}^{t''}Z_{t',t''}\subset
\rho_{t}^{t''}Z^{q-1}(t,C^{\hidot}(\sigma)) +
d
F^{q-2}(t'')$$
 L'application  $C' \to Z_{t', t''}$,
$(\xi,f) \mapsto ( r_{t'}^{t''})^{-1}
(\rho_t^{t''} \zeta +
df)$
est une contraction surjective. Elle poss\'ede une section
dont les
compos\'ees \`a gauche respectives $\Sigma_{t'}^{t}$
$\Phi_{t'}^{t''}$   avec  les deux applications
continues $C'\to
Z^{q-1}(C^{\hidot})(t)$ et $C' \to F^{q-2} (t'')$
v\'erifient $r_{t'}^{t''}=\rho_{t}^{t''}\Sigma_{t'}^{t} +
d\Phi_{t'}^{t''}$.

Les deux applications
${\cal A}$-lin\'eaires continues $S_{t'}^{t}= \Sigma_{t'}^{t}p, \
F_{t'}^{t''}=\Phi_{t'}^{t''}p$ v\'erifient les conditions requises.

\begin{lem} \label{mome}
Soient $0<t''<t'<t$ assez petits. Il existe un module
hilbertien de
type fini $M$, trois morphismes
de modules mont\'eliens $\mu: Z^{q-1}(C^{\hidot}(\sigma)) (t')\to
M$ , $\omega: M\to Z^{q-1}(C^{\hidot}(\sigma))(t')$ et $\tilde
F_{t'}^{t''}:Z^{q-1}(C^{\hidot}(\sigma)) (t')\to F^{q-2}(t'')$
tels que:
$$\rho_{t'}^{t''}|_{Z^{q-1}(t',C^{\hidot}(\sigma))} =
\rho_{t'}^{t''}\omega \mu +d\tilde F_{t'}^{t''}$$
\end{lem}
{{\bf Preuve: }}
 Le morphisme $\rho_t^{t'}$
 est ${\cal A}$-compact.
Par \ref{coroco}, il en va ainsi de
 $ c=\rho_{t}^{t'}S_{t'}^{t}$
qui est un endorphisme du module hilbertien
$Z^{q-1}(C^{\hidot}(\sigma))(t')$.
 Les points 3 et 4 du lemme \ref{lemco}
permettent de trouver un module hilbertien de type fini
$M$,
$\omega: M \to Z^{q-1}(C^{\hidot}(\sigma))(t')$
et $\tilde \mu:Z^{q-1}(C^{\hidot}(\sigma)) (t')\to M$
deux morphismes
tels que $c= \omega \tilde \mu +\delta$
o\`u  $\| \delta\|$
est un endomorphisme de
norme triple strictement inf\'erieure \`a $1$.
Soit $u =\sum_0^{+\infty} \delta^n$
l'inverse de $1-\delta$.

De l'\'equation
$\rho_{t''}^{t'} = \rho_{t''}^{t'} c +dF_{t'}^{t''}$,
on tire
$\rho_{t''}^{t'}(1-\delta) =
\rho_{t''}^{t'}\omega \tilde \mu + dF_{t'}^{t''}$ puis
$\rho_{t''}^{t'} = \rho_{t''}^{t'}\omega (\tilde \mu u) +
d(F_{t'}^{t''}u)$.$\qed$.

Fixons d\'esormais un triplet  $0<t''<t'<t$ assez petit et
appelons
$(\hat L^{\hidot},\hat \sigma^{\hidot})$ le complexe obtenu \`a partir
du couple $(M,\omega)$ du lemme \ref{mome} selon la
proc\'edure d\'ecrite ci dessous:

 $M$ est un module hilbertien de type fini muni d'un
morphisme
$\omega: M\to Z^{q-1} ( C^{\hidot}(\sigma))$. Composant
$\omega$ \`a gauche avec $Z^{q-1} (C^{\hidot}(\sigma)) \to
C^{q-1}(\sigma))$
puis avec les projections naturelles $C^{q-1}(\sigma))\to
F^{q-1}$
et $C^{q-1} (\sigma) \to L^q$ nous obtenons deux morphismes
$\hat d:M\to L^q$ et $\hat \sigma: M \to F^{q-1}$
v\'erifiant
$\sigma^q \hat d = d\hat \sigma$, i.e.: un complexe:

$$\hat L^{\hidot}:=\hat L^{\hidot} (M,\omega)= (\hat L^{q-1} :=M) 
\buildrel {\hat
d} \over
{\to}
L^{q} \buildrel{d}\over{\to} \ldots $$

et un morphisme de complexes $\hat \sigma^{\hidot}:=
\hat \sigma^{\hidot}(M,\omega): \hat L^{\hidot} \to
F^{\hidot}$
d\'efini par $\hat \sigma^i = \sigma ^i, i \ge q$ et $\hat
\sigma^{q-1} = -\hat \sigma$.

Observons que $C^i(\sigma)(s) = C^i(\hat \sigma)(s), \
i\ge q-1, s\le t'$.

\begin{lem}
$\rho_{*t}^{*t''} Z^{q-1}_*( t, C^{\hidot}(\sigma))
\subset \rho_{t'}^{*t''} Z^{q-1}( C^{\hidot}(\hat\sigma))( t')
\subset \partial C^{q-2}_*(\hat\sigma)(t'')$
\end{lem}

{\bf{Preuve :}} La premi\`ere inclusion r\'esulte de:
$$\rho_{*t}^{*t''} Z^{q-1}_*( t, C^{\hidot}(\sigma))
=\rho_{t'}^{*t''}\rho_{*t}^{t'}Z^{q-1}_*( t, C^{\hidot}(\sigma))
\subset \rho_{t'}^{*t''} Z^{q-1}(C^{\hidot}(\sigma))( t')  $$
 Ceci
joint au lemme \ref{mome} permet d'associer \`a tout
\'el\'ement $\zeta \in Z^{q-1}_*(t,C^{\hidot}(\sigma))$
un \'element $( \zeta', f) \in F^{q-2}(t'') \oplus M$
tel que $\rho_{*t}^{t''}\zeta= df+\omega m$.

Or, par d\'efinition, $C^{q-2}(\hat \sigma) = F^{q-2}
\oplus M$, $C^{q-1}(\hat \sigma) =C^{q-1} (\sigma)$ et la
diff\'erentielle de $C^{\hidot}(\hat \sigma)$ est donn\'e par la
formule $\partial_{\hat \sigma} (f,m)= \partial
_{\sigma}f+\omega m =df+\omega m$.

Par suite $\rho_{*t}^{*t''} \zeta =\partial_{\hat \sigma}
\rho_{t''}^{*t''}(f,m)\in \partial C^{q-2}(\hat
\sigma)_*(t'')$.$\qed$.

Soit $\theta <t''$.
\begin{eqnarray*}
Z_*^{q-1}(\theta,C^{\hidot}(\hat \sigma))&=&Z_*^{q-1}(\theta,
C^{\hidot}(\sigma))\\
&\subset & \rho_{*t}^{*\theta} Z^{q-1}_*(t,C^{\hidot}(\sigma))
+\partial C_*^{q-2}(\theta,\sigma)\\
& \subset& \rho^{*\theta}_{*t''}\partial C_*^{q-2}(t'',\hat
\sigma)+\partial C_*^{q-2}(\theta,\sigma)\\
& \subset& \partial C_*^{q-2}(\theta,\hat \sigma)
\end{eqnarray*}

La deuxi\`eme ligne r\'esulte du fait que l'identit\'e est
un qi-morphisme de $C^{\hidot}(\sigma)$. La derni\`ere
ligne fournit $H_*^{q-1}(\theta,C^{\hidot}(\hat \sigma))=0$.
Mais,
   $$i\ge q \Rightarrow H_*^i(\theta,C^{\hidot}(\hat \sigma))=0$$
puisque en ces degr\'es les cohomologies de $C^{\hidot}(\sigma)$
et $C^{\hidot}(\hat \sigma)$ sont les m\^emes. Utilisant \ref{qiq},
nous concluons que $\hat \sigma$ est un qi-morphisme
d'ordre $q$. Ceci conclut la preuve de la premi\`ere assertion
de la proposition \ref{tfini}.

La deuxi\`eme assertion se prouve par la m\^eme technique
\`a l'aide de la propri\'et\'e r\'ecursive suivante:

$P(q)$: Il existe $L^{\hidot}$ un
complexe de module hilbertiens et deux
morphismes de complexes $\psi^{\hidot}:T_{q-1} Q^{\hidot} \to L^{\hidot}$
et $\sigma^{\hidot}: L^{\hidot} \to F^{\hidot}$ un qi-morphisme d'ordre $q$
tels que $T_{q-1} \phi^{\hidot}= \sigma^{\hidot} \psi^{\hidot}$.

La preuve de $P(q) \Rightarrow P(q-1)$
 r\'esulte
des arguments pr\'ec\'edents.
Observons, en effet, que $\forall x
 \in Q^{q-1} \ d\psi^qdx=\psi^{q+1} ddx=0,\
d\phi^{q-1} x-\sigma^q\psi^qdx=
( \phi^q-\sigma^q\psi^q)dx=0$.
Par suite $e: Q^{q-1} \to C^{q-1}(\sigma)=
F^{q-1}\oplus L^q$ d\'efini par $e(x)=
 ( -\phi^{q-1}x,\psi^qdx)$ est \`a
 valeurs dans $Z^{q-1}(C^{\hidot}(\sigma))$.
Choisissons $t>t'>t''$ convenablement et observons que le
couple $(M,\omega)$ donn\'e par le lemme \ref{mome}
a \'et\'e construit comme une inclusion de sous-module ferm\'e:
$\omega: M\to Z^{q-1}(t',C^{\hidot}(\sigma))$. Nous d\'efinissons
un nouveau couple $(M',\omega')$ comme \'etant l'inclusion
du sous module ferm\'e
$\overline{M+eQ^{q-1}}$ qui est hilbertien
de type fini. Observons que l'inclusion $M'\subset M$
induit
un morphisme naturel de complexes
$\mu ^{\hidot}:\hat L^{\hidot}(M,\omega) \to \hat L^{\hidot}(M',\omega') $
tel que $\hat \sigma (M,\omega)
 =\hat \sigma (M',\omega')\mu^{\hidot}$. De cette
factorisation, on d\'eduit que $\hat \sigma (M',\omega')$
 est un qi-morphisme d'ordre $q-1$. Par ailleurs, on v\'erifie
ais\'ement que poser
$(\psi^{q-1}:Q^{q-1}\to \hat L ^{q-1}(M',\omega')):=
(e: Q^{q-1} \to M')$ d\'efinit un morphisme
$T_{q-1} Q^{\hidot} \to \hat L ^{\hidot}(M',\omega')$ avec
 les propri\'et\'es requises.$\qed$.

Soit $i\in \{b,-\}$.
Comme $E_f({\cal A})$
a assez de projectifs, les foncteurs naturels
$\mu^i:K^i(P_f({\cal A})) \to D^i(E_f({\cal A}))$  sont des \'equivalences 
de
cat\'egories.
D'autre part, nous avons d\'ej\`a implicitement rencontr\'e le
foncteur $\nu^i: K^i(P_f({\cal A}))\to K^i(MontMod_{\cal A})^{qhtf}$
et on note $\xi^i$ le foncteur de localisation:
 $$\xi^i: K^i(MontMod_{\cal A})^{qhtf}\to  K^i(MontMod_{\cal 
A})^{qhtf}_{Qi}$$

\begin{theo}\label{tfini2}
Le foncteur $\eta^i=\xi^i\nu^i(\mu^i)^{-1}$  r\'ealise une
\'equivalence de cat\'egories
$$\eta^i: D^i(E_f({\cal A}) )\to
K^i(MontMod_{\cal A})^{qhtf}_{Qi}$$
\end{theo}
{\bf{Preuve: }}Ce th\'eor\`eme  r\'esulte via \cite{Ha2}, prop. 3.3, p.33 
du premier point de la proposition  \ref{tfini}.

De \ref{tfini}, le lecteur ignorant tout des cat\'egories d\'eriv\'ees 
pourra d\'eduire le fait que,  pour tout  complexe de modules mont\'eliens 
qhtf born\'e sup\'erieurement $F^{\hidot}$, 
$\underrightarrow{\lim}_{t\to 0} H^q(t,F^{\hidot})$ s'obtient comme image
par le foncteur d'oubli $O$ 
d'un \'el\'ement bien  d\'etermin\'e de $E_f({\mathcal A})$. Moins 
pr\'ecis  que \ref{tfini2}, ce dernier  \'enonc\'e pourrait suffire.

\section{Th\'eor\`eme A de Cartan en Cohomologie $L^2$}

\subsection{R\'esolution de Dolbeault, acyclicit\'e locale des images 
directes $L^2$ et r\'esolution de \v Cech}

Soit $\tilde X$ une $\Gamma$-vari\'et\'e complexe.
Soit $\pi:\tilde X\to X= \Gamma \backslash \tilde X$.

Soit ${\cal V} \in V_{\bar \Gamma, \chi} (\tilde X)$. Soit
 $V$ le fibr\'e vectoriel associ\'e que
 l'on munit d'une m\'etrique hermitienne.

Soit $U$ un ouvert de $X$. On pose:
$${\cal D}_2^{0,q}({\cal V})(U) =\{
s\in L^2_{loc} (\pi^{-1} (U), V\otimes \Omega^{0,q})
: \ \forall K \subset\subset U \int_{\tilde K}
\|s\|^2+ \| \bar \partial s\|^2 < \infty \}
$$

Et on d\'efinit un complexe de faisceaux, que nous
appelerons
le complexe de Dolbeault $L^2$ local, par la formule:

$$ {\cal D}_2^{\hidot} ({\cal V}) =({\cal D}_2^{0,0}({\cal V})
\buildrel
{\bar \partial}\over{\to}
{\cal D}_2^{0,1}({\cal V})\to \ldots)$$

\begin{lem}\label{dolbeault}
$l^2\pi_* {\cal V} \to {\cal D}_2^{\hidot}({\cal V}) $
 est une r\'esolution
acyclique de $l^2\pi_* {\cal V}$.
\end{lem}
{\bf {Preuve:}} C'est clair.
Donnons pourtant les d\'etails.

Il existe une
 partition de l'unit\'e
lisse et $\Gamma$-invariante sur $\tilde X$
subordonn\'ee \`a l'image r\'eciproque
d'un recouvrement ouvert quelcoque de $X$.
Par suite,
 les faisceaux $ {\cal D}_2^{0,q}({\cal V}), \  q\ge 0$
sont mous et donc acycliques.

Une section $L^2_{loc}$ d'un fibr\'e holomorphe
 tu\'ee par l'op\'erateur
$\bar \partial$ est holomorphe,
 d'o\`u l'exactitude au premier terme.

L'exactitude au termes suivants r\'esulte
 des estim\'ees $L^2$ d'H\"ormander:

Soit $q>0$. Soit $x \in X$.
 Soit $\phi_x \in {\cal D}^{0,q}({\cal V})_x$.
On suppose $\bar \partial \phi_x=0$. On peut supposer que
$\phi_x$ est le germe en $x$ de $\phi$ une section $L^2$
de $V$ sur un ouvert $\Gamma$-Stein $\tilde U$.
$\tilde U =\Gamma U^0$ avec $U^0$ une composante connexe de stabilisateur 
fini $\Sigma$.
 Soit $P$ une fonction $\Gamma$-invariante,
$\Gamma$-exhaustive strictement psh et lisse sur $X$.
Soit $t> \min P$ et $x\in \tilde U_t= P^{-1}( ] -\infty, t] )$.
On peut s'arranger pour trouver une m\'etrique de K\"ahler
 sur $ U_t^0$ et compl\`ete
( $U_t^0$
est une  vari\'et\'e de Stein.)
 que, par moyennisation, on peut supposer
invariante par
le sous groupe fini $\Sigma\subset \Gamma$. Cette m\'etrique
se transporte  \`a une m\'etrique $\Gamma$-invariante
$\omega$
compl\`ete sur $\tilde U_t=\Gamma/ \Sigma \times U_t^0$.
Composant la fonction psh $P - \min_U P +1$ avec
une fonction
 convexe croissante
au comportement suffisamment explosif
on trouve une fonction plurisousharmonique $\psi$,
 $\Gamma$-invariante
telle
que $ \int_{\tilde U_t} e^{-\psi} \|\phi\|^2 \omega ^n
<+\infty$
et telle que ${\cal V} \otimes K_X^{-1}$ est
uniform\'ement strictement n\'egatif au sens de Nakano
pour la m\'etrique
hermitienne $h'=e^{-\psi} \|.\|^2_h (\omega ^n)$.

Utilisant la formule de
 Bochner-Kodaira-Nakano et la technique d'H\"ormander,
\cite{Dem} pp. 26-32, on trouve $\mu$ une $q-1$ forme
telle que $\bar \partial \mu =\phi$ et
$\int e^{-\psi} \|\mu \|^2_h (\omega ^n)< +\infty$
et en particulier $\mu \in L^2_{loc} (\tilde U_t)$.
Son germe en $x$, not\'e $\mu_x$,
v\'erifie  $\bar \partial \mu_x=\phi_x$.$\qed$.

\begin{coro}\label{thbvect}
Soit ${\cal V} \in V_{\bar \Gamma, \chi} ({\tilde
X})$. Si  $\tilde X$ est $\Gamma$-Stein,
alors $H^q_{(2)} (\tilde X,{\cal V})=0$ pour $q>0$.
\end{coro}
{\bf{Preuve:}} Par le lemme \ref{dolbeault}, on
peut utiliser le complexe de Dolbeault $L^2$
local
pour calculer ce groupe de cohomologie. Par ailleurs, la preuve de
\ref{dolbeault} fournit la conclusion d\'esir\'ee.
 $\qed$.

\begin{prop}\label{thbcoh0}
Soit ${\cal F} \in C_{\bar \Gamma, \chi} ({\tilde
X})$. Si  $\tilde X$ est $\Gamma$-Stein,
 $H^q_{(2)} (\tilde X,{\cal F})=0$ pour $q>0$.
\end{prop}
{\bf{Preuve :}}
Soit $\tilde V$ un ouvert $\Gamma$-Stein de $\tilde X$
tel que $\Gamma \backslash \tilde V \subset \subset \Gamma \backslash \tilde
X$.
 Soit $R^{\hidot}_{\cal F} \to {\cal F}$ une
r\'esolution p\'eriodique
localement libre finie  de ${\cal F}$ d\'efinie sur $\tilde V$. Par la
proposition
\ref{lpexact}, le complexe  de faisceaux suivant est exact:

$$\ldots \to l^2\pi_*R^{-1} _{\cal F} \to
l^2\pi_*R^0_{\cal F}
 \to l^2\pi_*{\cal F}\to 0$$

Par suite $H^q_{(2)}(\tilde V, {\cal F}) \simeq {\Bbb H}^q(V,
\ldots \to l^2\pi_*R^{-1}_{\cal F}  \to l^2\pi_*R^0_{\cal
F})$.

Le corollaire \ref{thbvect} implique que
$H^i(V,l^2\pi_*R^j_{\cal F})=0, i>0$.
La suite spectrale d'hypercohomologie (voir \cite{GH}, pp. 445-447) 
$(''E^{pq}_r)_r$
v\'erifie donc $''E^{pq}_1 = 0 $  si $q\not =0$  et de plus,
$''E^{p,0}_1=H^0(V,l^2\pi_* R^p_{\cal F})$.
De l\`a $''E^{p,0}_1 =0, \  p>0$. De ce fait
$(''E^{pq}_r)_r$ d\'eg\'en\'ere en $''E_2$ et
$''E_2^{p,q} =0 $ sauf \'eventuellement quand $p=0$
et $q\le 0$. Puis ${\Bbb H}^q(V,l^2\pi_*R^{\hidot})=0$
sauf \'eventuellement si $q\le 0$. En particulier
$H^q_{(2)}(\tilde V,{\cal F})=0, \ q>0$
 (et aussi bien s\^ur $''E_2^{p,q}=0,
\ (p,q) \not = (0,0)$, car un
faisceau n'a pas de cohomologie en degr\'es
n\'egatifs.)

Comme
$\underleftarrow{\lim}_{V\subset\subset U}
H^q_{(2)}(\tilde V,{\cal F})
 =H^q_{(2)}(\tilde U,{\cal F})$,
 le lemme est d\'emontr\'e.  $\qed$.

\begin{prop}\label{tbcoh}
Soit $\tilde Z$ un $\Gamma$-espace complexe $\Gamma$-Stein. Soit
${\cal F} \in C_{\bar \Gamma}(\tilde Z)$. Pour tout $q\ge 1$, 
$H^q_{(2)}(\tilde Z,{\cal
F})=0$.

\end{prop}
{\bf {Preuve:}} En effet, il existe un $\Gamma$-plongement ferm\'e
$i:\tilde Z \to \tilde X$ o\`u $\tilde X$ est une
$\Gamma$-vari\'et\'e $\Gamma$-Stein et il ressort de \ref{modif0},
\ref{modif}
que $H^q_{(2)}(\tilde Z,{\cal F}) =H^q_{(2)}(\tilde X,i_*{\cal F})$.
\ref{tbcoh} est donc un corollaire de \ref{thbcoh0}.

Soit ${\frak V}=\{\tilde V_{i} \}_{i \in A}$ un
recouvrement ouvert d'un espace topologique. On note
$PF(A)$ l'ensemble des parties finies de $A$.
Pour $\alpha=\{ i_1, \ldots, i_p\}\in PF(A)$,
 on note $V_{\alpha}= V_{i_1} \cap \ldots\cap V_{i_p}$.

\begin{coro}\label{cech}
Soit $\tilde Z$ un $\Gamma$-espace complexe et ${\cal F} \in C_{\bar 
\Gamma}(\tilde
Z)$.

Soit ${\frak U}=\{\tilde U_{i} \}_{i \in A}$ un
recouvrement ouvert $\Gamma$-invariant
 $\Gamma$-localement fini
\footnote{Ceci signifie, par d\'efinition, que le recouvrement de 
$Z=\Gamma\backslash\tilde Z$
donn\'e par $\{\Gamma\backslash\tilde U_{i} \}_{i \in A}$ }
de $\tilde Z$
par des ouverts $\Gamma$-Stein.

Le complexe de \v Cech dont les espaces de
cochaines sont d\'efinies par:

 $$C_2^p (\tilde {\frak U}, {\cal F}) =
 \oplus _{\alpha \in PF(A) |\alpha | =p+1}
 H^0 (\tilde U_{\alpha}, l^2\pi_*{\cal F})$$

et les diff\'erentielles par la formule usuelle
calcule
la cohomologie $L^2$ de ${\cal F}$.
\end{coro}

\subsection{ Structures topologiques sur $H^0_{(2)}(\tilde X,{\cal F})$}

\paragraph{Faisceaux localement libres sur un espace r\'eduit}
Soit $\tilde U$ un $\Gamma$-espace $\Gamma$-Stein et r\'eduit.
Soit $P$ une fonction strictement psh sur $\tilde U$
$\Gamma$-invariante et $\Gamma$-exhaustive.

On r\'eutilise  les notations du lemme
\ref{frechet} en posant $\tilde U_{t}= P^{-1} (]-\infty,t[)$.

 Soit ${\cal V}
 \in V_{\bar \Gamma,\chi} (\tilde U)$.
Fixons une m\'etrique hermitienne sur ${\cal V}$
d\'efinie sur un ouvert
de la forme $\tilde U_{t_0}, \ t_0 >0$.

On d\'efinit, pour tout $t<t_0$ ,
$L^2(\tilde U_t,{\cal V})$ comme l'espace des
sections holomorphes $L^2$ de $\cal V$ sur $\tilde U_t$
muni de la norme hilbertienne $\|.\|_t$
d\'efinie par l'int\'egration sur $\tilde U_t$
 du carr\'e de la norme d'une section.

\begin{lem}\label{modhilbv}
$(L^2(\tilde U_t,{\cal V}), \| . \|_t)$ est
un $W^*(\bar \Gamma, \chi)$-module hilbertien projectif
m\'etris\'e.
\end{lem}
{\bf{Preuve :}} Rappelons que si le groupe $G$
agit \`a gauche sur un ensemble $E$, l'action se
relevant \`a un fibr\'e vectoriel
 $F\to E$,
$G$ agit sur l'espace des sections $\phi$ de $E$
par $L_g \phi (x) = g \phi (g^{-1}x)$.

Cel\`a \'etant dit, \`a $s\in L^2(\tilde U_t,{\cal V})$, on
associe la fonction sur $\bar \Gamma$ \`a valeurs
dans l'espace de
 Hilbert $H_t=L^2(U_t^0,{\cal V})$:
$$ F_t(s) (\bar g) = L_{\bar g^{-1}} s | _{ U^0_t}$$

$F_t(s)$ est une fonction $L^2$ sur $\bar \Gamma$
\`a valeurs dans l'espace de Hilbert $H_t$.
L'application lin\'eaire r\'esultante
 $I_{t}:L^2(\tilde U_t,{\cal V}) \to L^2(\bar \Gamma)\hat \otimes H_t$
est un plongement isom\'etrique  $\bar \Gamma$-\'equivariant.

Toute fonction de la forme $F=F_t(s)$ v\'erifie
$\forall z \in S, \chi(z) F(\bar g z) = F(\bar g)$.
Par suite $I_t$ est \`a valeurs dans $E_{\chi} \otimes H_t$.

Ceci r\'ealise $L^2(\tilde U_t,{\cal V})$ comme
un $W^*(\bar \Gamma,\chi)$-module hilbertien projectif.

\begin{lem}
\label{montelv}$\forall t'<t$
 l'application de restriction
 $\rho_t^{t'} : L^2(\tilde U_t,{\cal V})\to L^2(\tilde U_{t'},{\cal V})$
est une contraction injective $W^*(\bar \Gamma,
\chi)$-compacte.
\end{lem}
{\bf {Preuve :}}
Pour $t'<t$ ,
on a un diagramme commutatif:
$$
\begin{array}{lcl}
L^2(\tilde U_t,{\cal V})& \buildrel{ \rho_{t}^{t'} }\over{\to}&L^2(\tilde 
U_{t'},{\cal V}) \\
\downarrow I_t & & \downarrow I_{t'} \\
E_{\chi}\otimes H_t &
\buildrel{Id_{E_{\chi}} \otimes r_{t}^{t'}} \over{\to}
& E_{\chi} \otimes H_{t'}
\end{array}
$$

$r_{t}^{t'}$ d\'esigne l'application de restriction
$H_t \to H_{t'}$. Par le th\'eor\`eme de Montel,
 $r_{t}^{t'}$ est compact (au sens
habituel, c'est \`a dire ${\Bbb C}$-compact
au sens de cet article). On d\'eduit  que
$\rho_{t}^{t'}$ est $W^*(\bar \Gamma,\chi)$-compact
en utilisant $|||Id \otimes a||| \le ||| a|||$.
$\rho_{t}^{t'}$ est injective gr\^ace
\`a l'unicit\'e du prolongement analytique.

\begin{coro} La donn\'ee $W=(\{L^2(\tilde U_t,{\cal V})\}_{0<t<t_0},\{ 
\rho_t^{t'} \}_{0<t'\le 
t<t_0})$
d\'efinit un ${W^*(\bar \Gamma,\chi)}$-module pr\'emont\'elien au sens de
la d\'efinition
\ref{modulana}.

De plus, pour tout $t_0>t>0$,
$W_*(t)=H_{(2)}^0(\tilde U_t,{\cal V})$.
\end{coro}

\paragraph{Faisceaux analytiques coh\'erents sur un espace r\'eduit}

\begin{lem}\label{montelcoh}
 Soit ${\cal F}\in C_{(\bar \Gamma, \chi)}(\tilde U)$.

$H_{(2)}^0 (\tilde U_t,{\cal F})$  poss\`ede une structure
canonique d'espace de Fr\'echet.

Soit $c$ une pr\'esentation localement libre
de ${\cal F}$
 sur $\tilde U$, c'est \`a dire un \'epimorphisme p\'eriodique de la forme
$ {\cal W}\buildrel{c}\over{ \to} {\cal F} \to
0$ o\`u $V \in V_{\bar \Gamma,\chi}(\tilde U)$ .

 La structure
d'espace de Fr\'echet de $H_{(2)}^0(\tilde U,{\cal F})$ peut \^etre 
d\'efinie par une famille
$(\|.\|^c_{t'})_{ t'<t}$
croissante de semi normes pr\'ehilbertiennes
$ \Gamma$-invariantes.

Le s\'epar\'e
compl\'et\'e de $(H_{(2)}^0(\tilde U_t,{\cal F}), \|. \|^c_{t'})$
est ind\'ependant de $t$,
 est un $(\bar \Gamma,\chi)$-module hilbertien projectif
s\'eparable  not\'e $L^2(\tilde U_{t'},{\cal F})_c$.

Les morphismes naturels  $\rho_{t'}^{t''}$
$L^2(\tilde U_{t'},{\cal F})_c\to L^2(\tilde U_{t''},{\cal F})_c$ induits 
par l'identit\'e de 
$H_{(2)}^0(\tilde U_t,{\cal F})$
sont continus et ${W^*(\bar \Gamma,\chi)}$-compacts.

La donn\'ee ${\cal L}^2 (\tilde U,{\cal F})_{c,P}=(\{L^2(\tilde U_t,{\cal 
F})_c\}_{0<t<t_0},\{ 
\rho_t^{t'} \}_{0<t'\le t<t_0})$
d\'efinit un ${W^*(\bar \Gamma,\chi)}$-module pr\'emontelien.

De plus, pour tout $t>0$ assez petit
$({\cal L}^2 (\tilde U,{\cal F})_{c,P})_*(t)
=H_{(2)}^0(\tilde U_t,{\cal F})$.
\end{lem}
{\bf{Preuve :}}
A cause du lemme \ref{montelv}, c'est le cas si
${\cal F}$ est localement libre, les semi-normes \'etant
de vraies normes.

Soit
${\cal K}$   le noyau de ${\cal W} \to {\cal F}$ et ${\cal V }\to {\cal K}$
un \'epimorphisme p\'eriodique, d\'efini sur $\tilde U_{t_0}$.

La proposition \ref{tbcoh}
pr\'esente $H_{(2)}^0 (\tilde U_t,{\cal F})$ comme le conoyau
du morphisme continu
$ H_{(2)}^0(\tilde U_t,{\cal V} ) \to H_{(2)}^0(\tilde U_t, {\cal W})$.

Gr\^ace \`a \ref{tbcoh} l' image
de ce morphisme est $H_{(2)}^0(\tilde U_t, {\cal K})$.
En vertu de
\ref{closedness} ce sous espace est ferm\'e pour la
topologie de Fr\'echet.
Le quotient d'un espace de Fr\'echet par un sous espace
ferm\'e est un espace de Fr\'echet. Ceci construit une structure de 
Fr\'echet sur 
$H_{(2)}^0(\tilde U_t, {\cal F})$.

 Posant pour $f\in H_{(2)}^0 (\tilde
U_t,{\cal F})$ :
$$\|f\|^c_{t'} =\inf _{w \in H_{(2)}^0 (\tilde U_t,{\cal W}), \ c(w)
=f}
\|w\|_{t'}$$
nous obtenons une seminorme $\|.\|^c_{t'}$ sur $H_{(2)}^0(\tilde
U_t,{\cal F})$. Elle est pr\'ehilbertienne, car elle
v\'erifie l'identit\'e du parall\'elogramme.

Consid\'erons l'espace vectoriel norm\'e
$(H_{(2)}^0(\tilde U_t,{\cal F})/ Ker\|.\|_{t'}^c,\|.\|_{t'}^c)$.
Son compl\'et\'e
est, par d\'efinition
$L^2(\tilde U_{t'},{\cal F})_{c,t}$.
$c$ d\'efinit une contraction continue et surjective
$$(H_{(2)}^0(\tilde U_t,{\cal V}), \|. \|_{t'})
\to
(H_{(2)}^0(\tilde U_t,{\cal F})/ Ker\|.\|_{t'}^c,\|.\|_{t'}^c)$$

Par suite $c$ se prolonge
aux compl\'et\'es de ces
 espaces vectoriels norm\'es
et d\'efinit
une surjection continue
 $L^2(\tilde U_{t'},{\cal V})
 \to L^2(\tilde U_{t'},{\cal F})_{c,t}$ qui munit
 $L^2(\tilde U_{t'},{\cal F})_{c,t}$
d'une structure de $(\bar \Gamma,\chi)$-module
 hilbertien gr\^ace \`a \ref{modhilbv} . Pour tous $t_1>t_2 >t'$
$H_{(2)}^0(\tilde U_{t_1},{\cal V})$ est dense dans
$H_{(2)}^0(\tilde U_{t_2},{\cal V})$ pour la norme $\|. \|_{t'}$,
ce qui implique que $
(L^2(\tilde U_{t'},{\cal F})_{c,t}, \|. \|_{t'}^c)$
est ind\'ependant de $t>t'$. On peut donc le noter
$(L^2(\tilde U_{t'},{\cal F})_{c}, \|. \|_{t'}^c)$.

Les autres points de l'\'enonc\'e sont
des cons\'equences directes de \ref{montelv},
\`a l'exception de
l'ind\'ependance de la pr\'esentation $c$ de la structure
de Fr\'echet sur $H_{(2)}^0(\tilde U_t,{\cal F})$. Il suffit de
comparer les structures de Fr\'echet associ\'ees \`a
deux pr\'esentations locales localement libres de ${\cal F}$
$c$ et $c'$ telles qu'il existe un morphisme de pr\'esentations $c'\to c$,
tous ces objets pouvant \^etre d\'efinis sur un voisinage de
$\tilde U_{t_0}$. Un morphisme de pr\'esentations est un diagramme de la 
forme:
$$
\begin{array}{cccccccc}
& & {\cal W}'& \buildrel{c'}\over{\to}& {\cal F} &\to &0 \\
& & \downarrow& & \| & & \\
& & {\cal W}& \buildrel{c}\over{\to}& {\cal F} &\to &0
\end{array}
$$
Dans ces conditions, pour tout $t'<t_0$, il existe $K>0$
 $\|. \|_{t'}^c\le K \|.\|_{t'}^{c'}$. Par suite:
$$Id: (H_{(2)}^0(\tilde U_t,{\cal F}),\{\|. \|^{c'}_{t'}\}_{t'>0})
\to
(H_{(2)}^0(\tilde U_t,{\cal F}),\{\|. \|^{c}_{t'}\}_{t'>0})$$
est continue et les deux structures de Fr\'echet
$F_{c'}$ et $F_{c}$ d\'efinies respectivement par $c$ et $c'$
coincident.  $\qed$.

Un point ennuyeux est que ceci ne signifie pas que
$\| .\|^c_{t'}$ est \'equivalente \`a $\|. \|^{c'}_{t'}$.
Ceci se traduit par le fait d\'esagr\'eable que le module mont\'elien
${\cal L}^2 (\tilde U,{\cal F})_{c,P}$ d\'epend  de $c$ et $P$.

\paragraph{Cas g\'en\'eral}

\begin{defi}\label{casgen}
Soit $\tilde Z$ un $\Gamma$-espace $\Gamma$-Stein,  $i:\tilde Z \to \tilde
X$ un plongement ferm\'e vers un $\Gamma$-espace $\Gamma$-Stein, $P$ une 
fonction
 strictement psh
lisse $\Gamma$-exhaustive sur $X$ et $c$ une pr\'esentation dans
$V_{\bar \Gamma}(\tilde X)$ de $i_*{\cal F}$.

On d\'efinit ${\cal L}^2(\tilde Z,i,P,c,{\cal F})$ comme le module
mont\'elien ${\cal L}^2(\tilde X,i_*{\cal F})_{c,P}$.
\end{defi}

\begin{lem}\label{lien}
${\cal L}^2(\tilde Z,i,P,c,{\cal F})_*(t) =H_{(2)}^0(\tilde Z_t,{\cal F})$
\end{lem}
{\bf{Preuve:}} Ceci r\'esulte de \ref{modif0}, \ref{modif} et
\ref{montelcoh}.

\paragraph{Structure de Fr\'echet}

\begin{prop}
Soit $\tilde X$ un $\Gamma$-espace complexe. Soit
${\cal F}$ un faisceau coh\'erent $(\bar \Gamma,
\chi)$-p\'eriodique sur $\tilde X$.

$H_{(2)}^0 (\tilde X,{\cal F})$  poss\`ede une structure
canonique d'espace de Fr\'echet. De plus, cette structure
d'espace de Fr\'echet peut \^etre d\'efinie par une famille
(non canonique) $(\|.\|_n)_{ n \in {\Bbb N}}$
croissante de pr\'enormes pr\'ehilbertiennes
$\bar \Gamma$-invariantes dont les s\'epar\'es
compl\'et\'es
sont des $(\bar \Gamma,\chi)$-modules hilbertiens
s\'eparables.
\end{prop}
{\bf{Preuve:}} Utilisant qu'un espace complexe est paracompact et
v\'erifie le second axiome de d\'enombrabilit\'e,
on munit les termes $C^0$ et $C^1$
du complexe de \v Cech de $l^2\pi_*{\cal F}$ associ\'e \`a un
recouvrement $\Gamma$-localement fini par ouverts
$\Gamma$-pseudoconvexe d'une telle structure de Fr\'echet.
La diff\'erentielle de \v Cech est alors continue et le
r\'esultat s'ensuit. Puisque cette proposition ne
servira pas dans la suite de ce texte,
la preuve
 de l'unicit\'e de la structure de Fr\'echet
ainsi fabriqu\'ee est laiss\'ee au lecteur.
\subsection{ Th\'eor\`eme de finitude}

Le but de ce paragraphe est d'associer au complexe de \v Cech $L_2$ de
\ref{cech}
un complexe de modules mont\'eliens qhtf qui lui soit quasi isomorphe. 
Nous n'avons pu mener \`a bien la 
construction sans une certaine quantit\'e de donn\'ees annexes 
qui rendent lourde la r\'edaction d'une construction ais\'ee dans son 
principe.

Soit $\tilde Z$ un espace complexe cocompact.

\paragraph{Constructions}

\begin{defi} Soit $Z$ un espace complexe. 
Soit
${\cal F}$ un faisceau analytique 
coh\'erent sur $ Z$.

Soit $i_l, \ l=1,2$ un plongement 
localement ferm\'e
vers l' espace complexe $Z_l$,
$p$ un morphisme tel que $i_1=p\circ i_2$,
$c_1 {\cal V}_1 \to i_{1*}{\cal F}$ une 
pr\'esentation.

La pr\'esentation $p^{-1}c_1:p^*{\cal V}\to 
i_{2*}{\cal F}$
est l'unique pr\'esentation telle que pour 
tous $z\in Z$,  $f\in {\cal F}_z, \ s\in {\cal
V}_{i_1(z)}$ tels que $c_1(s) =i_{1*}f$ on 
ait $p^{-1}c_1 (p^*s)=
i_{2*}f$.
\end{defi}

\begin{defi}
Soit $\tilde Z$ un $\Gamma$-espace 
complexe. $\hat C_{\bar \Gamma,\chi}(\tilde Z)$
 d\'esigne la cat\'egorie suivante:

\begin{itemize}
\item Les objets sont les quintuplets 
$({\cal F},\tilde U, i, P, c)$,
${\cal F}\in C_{\bar \Gamma,\chi}(\tilde Z)$, 
$\tilde U$ un ouvert
$\Gamma$-invariant de $\tilde Z$, $i:\tilde 
U\to \tilde X$ un plongement
 ferm\'e vers le $\Gamma$-espace complexe 
r\'eduit $ \tilde
X$, $P$ une $\Gamma$-fonction continue 
strictement psh et $\Gamma$-exhaustive
et $c:{\cal V }\to i_{*}{\cal F}$ une 
pr\'esentation
localement libre p\'eriodique.
\item Soient $X_l=({\cal F}_l,\tilde U_l, 
i_l, P_l, c_l)\  l=1,2$ deux
objets.

Il n'y a de  fl\'eche $X_1 \to X_2$ que
si $\tilde U_2\subset \tilde U_1$
auquel cas une fl\'eche
 est la donn\'ee d'un triplet
$(\phi,p,\lambda)$ o\`u $\phi:{\cal F}_1\to 
{\cal F}_2 $ est un $\bar \Gamma$-morphisme
$p:X_2\to X_1$ une application holomorphe 
telle que $i_1=p\circ
i_2$, $P_2\ge P_1\circ p$ et $\lambda: 
p^*{\cal V}_1 \to {\cal V}_2$
est un morphisme
p\'eriodique faisant commuter le diagramme:
$$
  \begin{array}{ccc}
    p^*{\cal V}_1 & \buildrel{p^{-1}c_1} 
\over{\longrightarrow}& i_{2*}{\cal F}_1\\
  \lambda  \downarrow & & 
i_{2*}\phi\downarrow \\
     {\cal V}_2 & \buildrel{c_2} 
\over{\longrightarrow}& i_{2*}{\cal
     F}_2
  \end{array}
$$\\
\item La composition des fl\'eches est 
donn\'ee par la
     formule \'evidente.
\end{itemize}
\end{defi}

Reformulons la d\'efinition \ref{casgen} et le lemme \ref{lien} comme suit:

\begin{lem} \label{cavecanem}
Il existe un foncteur covariant
 ${\cal L}^2: \hat C_{\bar 
\Gamma,\chi}(\tilde Z) \to 
MontMod_{W^*(\bar \Gamma,\chi)}$
 tel que pour toute fl\'eche $F$ de la forme
 $$
F: X_1=({\cal F}_1,\tilde U_1, i_1, P_1, c_1)
 \buildrel{(\phi,p,\lambda)}\over{\to} 
X_2=({\cal F}_2,\tilde U_2, i_2, P_2, c_2)
 $$
et tout $t>0$ assez petit,
 le morphisme induit par le lemme \ref{lien}
 $${\cal L}^2(F)_*(t): H_{(2)}^0( \{z \in 
\tilde U_1| P_1(z)<t \},{\cal F}_1)
 \to  H_{(2)}^0(\{z \in \tilde U_2| P_2(z)<t 
\},{\cal F}_2)$$
 soit
 \'egal \`a la composition du morphisme  
$H_{(2)}^0(\phi)$  et de la
 restriction de $U_1$ \`a $U_2$.
\end{lem}

\begin{defi} 
D\'efinissons une cat\'egorie $\wave C_{\bar\Gamma,\chi}(\tilde Z)$
par les donn\'ees suivantes:

\begin{itemize}
 \item Objets: Un objet de $\wave C_{\bar \Gamma,\chi}({\tilde Z})$
est une famille d'objets de $\hat C_{\bar \Gamma}({\tilde Z})$,
 $D=(X_l=({\cal F}_l,\tilde U_l, i_l, P_l, 
c_l))_{l\in A}$
 sujette aux restrictions suivantes:
 \begin{itemize}
\item $X_l\in\hat C_{\bar \Gamma,\chi} (\tilde Z)$
\item Il existe un faisceau analytique coh\'erent p\'eriodique ${\cal 
F}={\cal F}_D\in C_{\bar \Gamma,\chi}(\tilde 
Z)$, tel que $\forall l \ \ \ {\cal F}_l={\cal F}$.
\item $(\tilde U_l)_l$ est un recouvrement 
$\Gamma$-localement fini de $\tilde Z$
\item $(\{z\in \tilde U_l | P_l(z)<0 \})_l$ 
est aussi un
recouvrement de $\tilde Z$.
 \end{itemize}
\item Morphismes: Une fl\'eche $F:(X^1_l)_{l\in A} \to (X^2_n)_{n\in B}$
est la donn\'ee d'un morphisme p\'eriodique
$\phi:{\cal F}_{D^1} \to {\cal F}_{D^2}$, d'une application $\rho: B\to A$ 
et de $\hat C_{\bar\Gamma,\chi}(\tilde Z)$-fl\'eches $F_n:X^1_{\rho(n)}\to 
X^2_{n}$, avec $F_n=(\phi,p_n,c_n)$ et $\rho$ d\'efinit un raffinement 
$(\{z\in \tilde U^2_l | P^2_l(z)<0 \})_{l\in B}
\to (\{z\in \tilde U^1_l | P^1_l(z)<0 \})_{l\in A}$.
\item La composition de deux fl\'eches est d\'efinie par la formule 
\'evidente. 
\end{itemize}
\end{defi}

Soit $D$ un objet de $\wave C_{\bar \Gamma,\chi}({\tilde Z})$. 

Pour $\alpha \in PF(A)$ on d\'efinit un 
objet $X_{\alpha}
=({\cal F}_{\alpha},\tilde U_{\alpha}, 
i_{\alpha}, P_{\alpha}, c_{\alpha})$
de $\hat C_{\bar \Gamma}(\tilde Z)$ en 
posant:
\begin{itemize}
\item ${\cal F} _{\alpha}={\cal F}_D$
\item $\tilde U_{\alpha}= \cap_{l\in 
\alpha} \tilde U_l$
\item $i_{\alpha} = \times_{l\in \alpha} i_l$
\item $P_{\alpha} = \max _{l\in\alpha}P_l$
\item $c_{\alpha} : {\cal V}_{\alpha} \to 
i_{\alpha *} {\cal F} =
 \sum_{l\in\alpha} p^{-1}_{l,\alpha}c_l:
 \oplus_{l\in\alpha} p^*_{l,\alpha}{\cal 
V}_l \to p_{i,\alpha*}{\cal F}_D$
 o\`u $p_{l,\alpha}: \times_{k\in \alpha} 
X_k \to X_l$ est la
 projection naturelle.
\end{itemize}

Par construction, pour tout 
$(\alpha,\alpha')\in PF(A)$ tel que
$\alpha'\subset \alpha $, le triplet 
$\rho_{\alpha} ^{\alpha'}=(\phi_{\alpha} 
^{\alpha'},
 p_{\alpha} ^{\alpha'}, \lambda_{\alpha} 
^{\alpha'})$ avec
 $\phi_{\alpha} ^{\alpha'}=Id$
 $p_{\alpha} ^{\alpha'}= \times 
_{l\in\alpha'} p_{l,\alpha}$
 $\lambda_{\alpha} ^{\alpha'}  $ est 
l'inclusion naturelle
 $\sum_{\l\in \alpha'} p_{l,\alpha}^*{\cal 
V}_l\to \sum_{\l\in \alpha} 
p_{l,\alpha}^*{\cal V}_l$
est un $\hat C_{\bar \Gamma}(X)$-morphisme.

De plus, on a la condition 
$\rho_{\alpha'}^{\alpha''}\rho_{\alpha}^{\alpha'}
=\rho_{\alpha}^{\alpha''}$.

On d\'efinit un ${W^*(\bar 
\Gamma,\chi)}$-module mont\'elien
gr\^ace \`a \ref{montelcoh}, en posant:
 $${\cal L}^2 C^p(D)=
\oplus _{|\alpha|=p+1} {\cal 
L}^2(X_{\alpha})$$

Imitons la d\'efinition de la 
diff\'erentielle de \v Cech en
d\'efinissant $\delta:{\cal L}^2 C^p 
(D)\to {\cal L}^2 C^{p+1}(D)$ comme
le morphisme dont la matrice est
$(\epsilon_{\alpha,\alpha'} {\cal 
L}^2(\rho_{\alpha}^{\alpha'}))_{\alpha,\alpha'
}$, avec 
$\epsilon_{\alpha,\alpha'}=\pm 1$, le signe \'etant donn\'e par la r\'egle 
usuelle pour la diff\'erentielle de \v Cech.

Par \ref{montelcoh}, $({\cal L}^2C^{\hidot}(D),\delta)$ est un complexe 
born\'e de
modules mont\'eliens.

\begin{lem}\label{435}
$({\cal L}^2 C^{\hidot}(D),\delta)$
est un complexe de modules mont\'eliens qhtf 
et pour tout $t>0$ assez petit
$H^{\hidot}_*(t,({\cal L}^2 C^{\hidot}(D),\delta))=H^{\hidot}_{(2)}(\tilde 
Z,{\cal F}_D)$.
\end{lem}
{\bf{Preuve :}} En effet,  gr\^ace \`a 
\ref{lien},  $({\cal L}^2 
C^{\hidot}(D)_*(t),\delta)$
s'identifie au complexe de \v Cech $L_2$ de ${\mathcal F}$ 
relatif au recouvrement
$(\{z\in U_l| P_l(z)<t\}_l )$ (cf. la d\'efinition du corollaire 
\ref{cech}) .

Par construction, on a:
\begin{lem}
L'assignement $D\mapsto ({\cal L}^2 C^{\hidot}(D),\delta)$ est sous jacent 
\`a un 
foncteur covariant $\wave C _{\bar \Gamma,\chi}(\tilde Z) \to 
C^b(MontMod_{W^*(\bar 
\Gamma,\chi)})^{qhtf}$. 
\end{lem}

\paragraph{Unicit\'e, Existence, Fonctorialit\'e}

Ce paragraphe utilise les notations du th\'eor\`eme \ref{tfini2}.

\begin{lem}\label{441}
Soit ${\cal F}$ un objet de $C_{\bar \Gamma,\chi}(\tilde Z)$.
Pour tout couple $D^1
,D^2$ d'objets de $\wave C_{\bar\Gamma,\chi}(\tilde Z)$,
tel que ${\cal F}_{D^i}= {\cal F}$, on 
peut
construire un isomorphisme entre
$({\cal L}^2 C^{\hidot}(D^1),\delta)_{Qi}$ et
 $({\cal L}^2 C^{\hidot}(D^2),\delta)_{Qi}$, {\em 
unique} dans
  $K^b(MontMod_{W^*(\bar 
\Gamma,\chi)})^{qhtf}_{Qi}$.
\end{lem}
{\bf{ Preuve:}} En effet, on peut toujours trouver un objet $D^3$
tel que ${\cal F}_{D^3}={\cal F}$
et des $\wave C_{\bar\Gamma,\chi}(\tilde Z)$-fl\'eches $D^i\to D^3$
dont le $C_{\bar \Gamma,\chi}(\tilde Z)$ morphisme associ\'e est 
l'identit\'e. 
La fl\'eche $({\cal L}^2 C^{\hidot}(D_1),\delta)\to ({\cal L}^2 
C^{\hidot}(D_3),\delta)$
est un qi-morphisme puisqu'elle induit des isomorphismes sur les groupes 
de  cohomologie 
par \ref{435}. Ce qi-morphisme d\'efinit apr\'es localisation un 
isomorphisme en cohomologie qui est ind\'ependant de la fl\'eche $D_1\to 
D_3$ choisie
puisqu'il rel\`eve les applications de raffinement entre 
 complexes de \v Cech $L_2$ de ${\cal F}$. $\qed$

Soit $\cal F$ un faisceau coh\'erent $(\bar\Gamma,\chi)$-p\'eriodique. 
Utilisant \ref{condition3}, on construit ais\'ement un objet $D({\cal F})$ 
de 
$\wave C_{\bar \Gamma}(\tilde Z)$ avec ${\cal F}_{D({\cal F})}={\cal F}$. 
Toujours avec \ref{condition3}, il est possible de relever toute
fl\'eche ${\cal F} \to {\cal G}$ en une fl\'eche
$D({\cal F}) \to D({\cal G})$. En raison de \ref{441}, 
on prouve le:

\begin{theo} \label{tfini3}
Soit $\tilde X$ un 
$\Gamma$-espace complexe
cocompact.

Il existe un foncteur triangul\'e (unique 
\`a \'equivalence pr\`es)

 $$R:D^i(C_{\bar\Gamma,\chi}(\tilde X)) \to
D^i(E_f(W^*(\bar\Gamma,\chi)) \ i=b,-$$

 sp\'ecialisant au $\delta$-foncteur
$( H_{(2)} ^q(\tilde X, .),\delta)_{q\ge 0}$
de \ref{etvoila}
 par oubli de structure (cf. \ref{isomalg}). Le foncteur $R$ est le 
foncteur d\'eriv\'e du foncteur sections globales $L_2$. 
\end{theo}

Soit $H^q_2(\tilde X,{\cal F})$  l'objet de $E_f(W^*(\bar\Gamma,\chi))$ 
d\'efini par $H^q_2(\tilde X,{\cal F})=H^q( R({\cal F}))$. 
L'assignement ${\cal F} \mapsto H^q_2(\tilde X,{\cal F})$ 
est sous jacent \`a un
 $\delta$-foncteur $(H^q_{2}(\tilde X, .),\delta)$ 
d\'efini sur $C_{\bar\Gamma,\chi}(\tilde Z)$
et \`a valeurs dans $E_f({\Gamma})$
 relevant \`a le $\delta$-foncteur \ref{defonct}
$(H^q_{2}(\tilde X, .),\delta)$, qui lui est \`a valeurs  \`a valeurs dans 
$E_f(W^*(\bar\Gamma,\chi))$

Si $\tilde X$ n'est pas suppos\'e cocompact, le th\'eor\`eme pr\'ec\'edent 
est valable dans le cas des faisceaux de support $\Gamma$-compact. Adapter 
la preuve donn\'ee ci dessus est trivial.

Pour justifier
 le titre de cet article, donnons la 
d\'efinition:

\begin{defi}
Soit $(\tilde X, \Gamma,\rho)$ un 
$\Gamma$-espace complexe.

Les invariants de Von Neumann
 d'un faisceau analytique coh\'erent
$(\bar \Gamma,\chi)$-p\'eriodique ${\cal 
F}$ sont:
\begin{itemize}
\item Les nombres de Betti $L^2$ de ${\cal 
F}$,
 $h^q_2 (\tilde X,{\cal F})
:= \dim_{W^*_l(\bar \Gamma, \chi)}
 H^q_{2} (\tilde X,{\cal F})\in {\Bbb R}^+$
\item Les invariants de Novikov-Shubin de 
${\cal F}$,
$NS^q(\tilde X,{\cal F}):=
NS(H^q_{2} (\tilde X,{\cal F}))\in NS_d$
\end{itemize}
\end{defi}

\section{Th\'eor\`eme d'indice $L_2$ d' Atiyah}

Cette section utilise  certaines id\'ees de la preuve de Forster-Knorr du 
th\'eor\`eme de coh\'erence de Grauert, voir \cite{GR} pp.188-207 pour 
justifier une suite spectrale de Leray-Serre en cohomologie $L_2$, qui est 
l'outil principal pour la formule de Riemann-Roch.

\subsection{Suite spectrale de Leray}
\paragraph{Syst\`emes de faisceaux relatifs \`a un recouvrement}

Soit $X$ un espace topologique 
et ${\frak U}=(U_i)_{i\in A}$ un recouvrement ouvert fini de $X$. On peut 
construire la cat\'egorie ab\'elienne
$SMod({\frak U})$ suivante:

Les objets
sont des donn\'ees de la forme $({\cal 
F}_{\alpha},\rho_{\alpha'}^{\alpha})_{\alpha, \alpha'\in PF(A)}$, tels que
\begin{itemize}
\item
${\cal F}_{\alpha}$ est un faisceau en groupes  ab\'eliens sur $U_{\alpha}$
\item
  $\rho_{\alpha}^{\alpha'}:{\cal F}_{\alpha}|_{\tilde U_{\alpha'}}\to 
{\cal F}_{\alpha'}$ est un morphisme d\'efini d\`es que  $\alpha \subset 
\alpha'$
\item   Ces donn\'ees v\'erifient
la condition de cocycle 
$\rho_{\alpha}^{\alpha'}\rho_{\alpha'}^{\alpha''}=\rho_{\alpha}^{\alpha''}$
\end{itemize}
Un  morphisme $\phi:  ({\cal 
F}_{\alpha},\rho_{\alpha'}^{\alpha})\to ({\cal 
G}_{\alpha},\rho_{\alpha'}^{\alpha})$ est une collection 
$(\phi_{\alpha})_{\alpha \in PF(A)}$
de morphismes de faisceaux $ \phi_{\alpha}: {\cal F}_{\alpha}\to {\mathcal 
G}_{\alpha}$ commutant aux morphismes de transition. 
Les r\'egles de composition sont claires. 

Il y a un foncteur naturel $n:Mod(X) \to SMod({\mathfrak U})$.
Au niveau des objets, il associe au  faisceau $F$ le syst\`eme de faisceaux 
$(F|_{U_{\alpha}},r_{\alpha}^{\alpha'})$,
$r_{\alpha}^{\alpha'}$ \'etant l'application de restriction. Ce foncteur 
identifie $Mod(X)$ \`a une sous cat\'egorie pleine de $SMod({\frak U})$.

Soit $U$ un ouvert de $X$, on note par $I^U$ l'inclusion de $U$ dans $X$.
Soit $\Phi=({\cal 
F}_{\alpha},\rho_{\alpha'}^{\alpha})_{\alpha, \alpha'\in PF(A)}$ un objet 
de $SMod({\mathfrak U})$. On d\'efinit un complexe de faisceaux sur 
$\tilde X$, en d\'efinissant les cochaines par la formule:
  $$K^n({\frak U}, F)= \oplus _{\alpha \in PF(A),  |\alpha|=n+1}
i_*^{U_{\alpha}} F_{\alpha}$$
et la diff\'erentielle par la formule de \v Cech.

Si $F$ est un faisceau ab\'elien, tel qu'il existe une base d'ouverts $U$
 pour lesquels la restriction de  $F$ \`a $U\cap U_{\alpha}$ est un 
faisceau acyclique d\`es que $\alpha\in PF(A)$ est une partie non vide, le 
morphisme naturel 
$F\to K^{\hidot}({\mathfrak U}, n(F))$ est un quasi isomorphisme.

Plus g\'en\'eralement, si $ F^{\hidot}$ est un complexe de tels faisceaux, 
on a aussi un quasi isomorphisme $F^{\hidot} \to 
K^{\hidot}({\frak U}, n(F^{\hidot}))$ \footnote{
Pr\'ecisons  que $K^{\hidot}({\frak U}, n(F^{\hidot}))$
est le complexe associ\'e au bicomplexe que la construction
pr\'ec\'edente fabrique. }.

Si $A\subset Mod(X)$ est une sous cat\'egorie de la cat\'egorie des 
faisceaux
en groupes ab\'eliens on peut d\'efinir la cat\'egorie $SA$ construite 
comme ci dessus en demandant que objets et morphismes 
       soient dans $A$. Par exemple si $A=C_{\bar \Gamma,\chi}(\tilde X)$, 
on a la cat\'egorie $SC_{\bar\Gamma,\chi}({\mathfrak U})$.

\paragraph{ Extension de $l^2\pi_*$ \`a 
une cat\'egorie d\'eriv\'ee convenable}

Soit $\tilde X$ un $\Gamma$-espace complexe
 et soit $\pi: \tilde X \to X$. Nous supposerons de plus que $X$
poss\'ede un recouvrement fini $\tilde{\frak U} $ par des
ouverts  $\Gamma$-Stein assez petits
\footnote{ Un ouvert $\Gamma$-Stein $\tilde U$ est {\em assez petit} s'il 
existe un autre ouvert $\Gamma$-Stein
$\tilde V$ tel que $\Gamma \backslash \tilde U \subset \subset 
\Gamma \backslash \tilde V$.}.

$SC_{\bar\Gamma,\chi}({\mathfrak U})$ contient une sous cat\'egorie pleine
$S$ naturellement \'equivalente \`a  $C_{\bar\Gamma,\chi} (\tilde X)$.

Nous allons \'etendre le foncteur $l^2\pi_*$ 
\`a $C^i (SC_{\bar\Gamma, \chi}({\mathfrak U}))$.

 Soit $\tilde j: \tilde V \subset 
\tilde X$
un plongement ouvert. Soit  ${\cal F}\in C_{\bar \Gamma,\chi} (\tilde V)$.
Consid\'erons  le faisceau 
$\tilde j_*{\cal F}$. On peut alors poser:

$$l^2\pi_* \tilde j_*{\cal F}= j_* l^2\pi_*{\cal F} $$

 \`A ${\cal F}^{\hidot}$
un objet de $C^i(SC_{\bar\Gamma, \chi} 
({\mathfrak U}))$, on  associe par ce proc\'ed\'e  le complexe
 $l^2\pi_* K^{\hidot} ({\mathfrak U},{\cal F}^{\hidot} )$, qui est un 
objet de $C^i(Mod_{W^*(\bar\Gamma,\chi)}(X))$. Cet assignement  est 
sous-jacent \`a un 
foncteur. On a, de plus, la:

\begin{prop}
Le foncteur $l^2\pi_*: C_{\bar\Gamma,\chi} (\tilde X)\to 
Mod_{W^*(\bar\Gamma,\chi)}(X)$ se prolonge \`a un foncteur exact
$D^iSC_{\bar\Gamma,\chi}({\mathfrak U}) \to 
D^i Mod_{W^*(\bar\Gamma,\chi)}(X)$. 
\end{prop}

La preuve du th\'eor\`eme \ref{tfini} implique aussi la:

\begin{prop}

Le foncteur $R\Gamma 
\circ l^2\pi_*$
s'\'etend \`a un foncteur triangul\'e $D^i_{C_{\bar \Gamma,\chi}(\tilde 
X)}  
S_{\bar\Gamma,\chi}({\mathfrak U}) \to D^iE_f(W^*(\bar\Gamma,\chi))$

\end{prop}

\paragraph{Th\'eor\`eme de l'image directe de Grauert}

Il r\'esulte de la preuve de Forster-Knorr du th\'eor\`eme de l'image 
directe de Grauert
(cf. \cite{GR} ch. 10.4), 
que, $\tilde f: \tilde Y \to \tilde X$ 
d\'esignant un morphisme propre $\Gamma$-\'equivariant ,${\cal 
G}^{\hidot}$ un complexe born\'e de $C_{\bar\Gamma,\chi}(\tilde Y)$,
et ${\mathfrak V}$ un recouvrement fini de $\tilde Y$ par des ouverts 
$\Gamma$-Stein assez petits, le complexe de \v Cech relatif
$\tilde f_* { C}^{\hidot}({\frak V},  {\cal G}^{\hidot})$ a la 
propri\'et\'e G suivante:

\begin{defi}
Soit $F^{\hidot}$ un complexe born\'e de $Mod_{\bar\Gamma,\chi}(\tilde X)$.
On dit que $F^{\hidot}$ a la propri\'et\'e G si,
pour $\tilde i^U: \tilde U \subset \tilde X$ un ouvert $\Gamma$-Stein 
assez petit, on peut  trouver un complexe born\'e 
 ${\cal F}^{\hidot}$ dans $C_{\bar \Gamma,\chi}(\tilde U)$ et un morphisme 
de complexes 
${\cal F}^{\hidot} \to F^{\hidot}|_U $
induisant un isomorphisme en cohomologie. 
\end{defi} 

Soit $F^{\hidot}$ un complexe de faisceaux v\'erifiant la propri\'et\'e G.
La technique de \cite{GR}, 10.4 pp.202-204, permet de 
construire un objet de $C^i(SC_{\bar\Gamma,\chi}({\mathfrak U}))$ 
$({\cal F}^{\hidot}_{\alpha}, \rho_{\alpha}^{\alpha'})$ 
et pour 
$ \alpha\in PF(A)$ des morphismes de complexes, induisant
des isomorphismes en cohomologie et faisant commuter
le diagramme:
$$
\begin{array}{ccc}
{\cal F}_{\alpha'}^{\hidot}|_{\tilde U_{\alpha}}&\to&  {\cal 
F}^{\hidot}_{\alpha} \\
 \downarrow& &\downarrow \\
F^{\hidot}|_{\tilde U_{\alpha}}& \buildrel{id}\over{\to}& 
F^{\hidot}|_{\tilde U_{\alpha}}
\end{array}
$$

Cette donn\'ee est en fait  un quasi isomorphisme dans 
$C^i(SMod_{\bar\Gamma,\chi}({\mathfrak U}))$:
$
{\cal F}^{\hidot} \to n(F^{\hidot})
$.

Ceci donne lieu \`a  un quasi isomorphisme
$$K^{\hidot}({\mathfrak U}, {\cal F}^{\hidot})\to
K^{\hidot}({\mathfrak U},{ F^{\hidot}})$$

Appliquant ces consid\'erations \`a 
$F^{\hidot}=\tilde f_* { C}^{\hidot}({\frak V},  {\cal G}^{\hidot})$, on 
prouve:

\begin{lem} Soit  $i \in \{ b,-\}$. Le foncteur d\'eriv\'e:
$$R\tilde f_*: D^i C_{\bar\Gamma,\chi}(\tilde Y)\to 
D^i_{C_{\bar\Gamma,\chi}(\tilde X)}Mod_{\bar\Gamma,\chi}(\tilde X)$$
factorise
 \`a un foncteur triangul\'e 
$$R\tilde f_*: D^i C_{\bar\Gamma,\chi}(\tilde Y)\to 
D^i_{C_{\bar\Gamma,\chi}(\tilde X)} SC_{\bar\Gamma,\chi}({\mathfrak U})$$

\end{lem}

Un int\'eret de ce lemme est qu'il permet  de d\'efinir
$l^2\pi_* R\tilde f_* {\cal G}^{\hidot}$, ${\cal G}^{\hidot}$
\'etant un objet de $D^i C_{\bar\Gamma,\chi} (\tilde Y)$. 
Il permet aussi de formuler la:

\begin{prop}\label{tautologie} Soit $\tilde f: \tilde Y \to \tilde X$ un 
morphisme 
propre \'equivariant, alors
$R f_* \circ l^2\pi_* = l^2\pi_* \circ R\tilde f_*$.
\end{prop}

{\bf Preuve :} 
Soit ${\cal G}^{\hidot}$ un complexe born\'e (resp. born\'e 
sup\'erieurement) de $C_{\bar \Gamma,\chi}(\tilde Y)$.
Soit un quasi isomorphisme 
${\cal F}^{\hidot} \to  n( f_* C^{\hidot}({\mathfrak V},{\cal 
G}^{\hidot}))$ avec ${\cal F}^{\hidot}$ un objet de 
$SC_{\bar\Gamma,\chi}({\mathfrak U})$.

On a le diagramme suivant de morphismes de complexes de faisceaux:

$$
\begin{array}{cccccccccc}
l^2\pi_*K^{\hidot}({\mathfrak U},{\cal F}^{\hidot}) &   
 & & & 
K^{\hidot}(\pi({\mathfrak U}), f_*C^{\hidot} ({\mathfrak V},l^2\pi_* {\cal 
G}^{\hidot}))& 
 & & & &\\
i_1\downarrow& & & &i_2 \downarrow & & & & & \\
\pi_*K^{\hidot}({\mathfrak U},{\cal F}^{\hidot})    &\buildrel 
{q}\over{\to} &\pi_*K^{\hidot}({\mathfrak U},f_* C^{\hidot}({\mathfrak 
V},{\cal G}^{\hidot})) & =& 
K^{\hidot}(\pi({\mathfrak U}), f_*C^{\hidot} ({\mathfrak V},\pi_* {\cal 
G}^{\hidot}))& & & & & 
\end{array}
$$

$i_1, i_2$ sont des inclusions et $q$ un quasi isomorphisme. Comme:
\begin{eqnarray*}
l^2\pi_*R\tilde f _*{\cal G}^{\hidot} & =&l^2\pi_*K^{\hidot}({\mathfrak 
U},{\cal F}^{\hidot}) \\
 Rf_* l^2\pi_*{\cal G}^{\hidot}& =& K^{\hidot}(\pi({\mathfrak U}), 
f_*C^{\hidot} ({\mathfrak V},l^2\pi_* {\cal G}^{\hidot}))
\end{eqnarray*}

\ref{tautologie} r\'esulte du lemme:

\begin{lem}
L'image de $q\circ i_1$ est contenue dans celle de $i_2$ et
le morphisme de complexes de faisceaux r\'esultant
$$l^2\pi_*K^{\hidot}({\mathfrak U},{\cal F}^{\hidot})
\to K^{\hidot}(\pi({\mathfrak U}), f_*C^{\hidot} ({\mathfrak V},l^2\pi_* 
{\cal G}^{\hidot}))$$
est un quasi isomorphisme.
\end{lem}

{\bf{Preuve :}} Cette question est de nature locale. Par suite, il suffit 
de traiter le cas o\`u $\tilde X$ est $\Gamma$-Stein et o\`u il existe 
un quasi isomorphisme global ${\cal F}^{\hidot}\to \tilde f_*C^{\hidot}
({\mathfrak V}, {\cal g}^{\hidot})$. On est alors ramen\'e au cas o\`u 
${\frak U}$ est
le recouvrement trivial $\{ \tilde X \}$. Red\'ecrivons la situation:

On a des inclusions (au sens naif) 
$l^2\pi_*{\cal F}^{\hidot} \to \pi_*{\cal F}$, 
$\pi_*{\cal F}^{\hidot} \to \pi_*\tilde f_* C^{\hidot}({\mathfrak V},{\cal 
G}^{\hidot})$. 

On a \'egalement une \'egalit\'e:
$$\pi_*\tilde f_* C^{\hidot}({\mathfrak V},{\cal G}^{\hidot}) =f_*\pi_* 
C^{\hidot}({\mathfrak V},{\cal G}^{\hidot})=
f_* C^{\hidot}({\mathfrak V},\pi_*{\cal G}^{\hidot})$$  

D'o\`u deux inclusions:
 
$$
l^2\pi_*{\cal F}^{\hidot} \to f_* C^{\hidot}({\mathfrak V},\pi_*{\cal 
G}^{\hidot}) \leftarrow f_* C^{\hidot}({\mathfrak V},l^2\pi_*{\cal 
G}^{\hidot})
$$

Il s'agit de montrer qu'en fait:
$$
l^2\pi_*{\cal F}^{\hidot}\subset f_* C^{\hidot}({\mathfrak 
V},l^2\pi_*{\cal G}^{\hidot})
$$

et que cette inclusion est un quasi isomorphisme.

 Le probl\`eme 
\'etant de nature locale il sera suffisant de montrer que, apr\`es 
application du foncteur des sections globales sur $X=\Gamma\backslash 
\tilde X$
$$
A^{\hidot}_2=\Gamma (  X,l^2\pi_*{\cal F}^{\hidot})\subset B^{\hidot}_2=
\Gamma (  X, f_* C^{\hidot}({\mathfrak 
V},l^2\pi_*{\cal G}^{\hidot}))
$$
et que cette inclusion est un quasi isomorphisme.

Posons $\tilde X =\Gamma X^0$ avec $X^0$ une  r\'eunion fine de 
composantes connexes Stein de groupe
d'isotropie $\Sigma$ fini. Posons  $Y^0=\tilde f^{-1} (X^0)$ et 
d\'esignons par  ${\frak V}^0$ le recouvrement induit par ${\mathfrak V}$.

Il y a un morphisme de complexes en espaces de Fr\'echet: 

$$A^{\hidot}=\Gamma(X^0,{\cal F}^{\hidot})
\to B^{\hidot}=\Gamma(X^0, \tilde f_*C^{\hidot}({\frak V}^0,{\cal 
G}^{\hidot}))$$

Ce morphisme est un quasi isomorphisme puisque les deux complexes ont
pour cohomologie $\Gamma(X^0,R^{\hidot}\tilde f _*{\cal G}^{\hidot})$.
De plus, ces complexes en espaces de Fr\'echet sont stricts \footnote{ 
C'est \`a dire tels que tels que l'image de
la diff\'erentielle est ferm\'ee.} et la structure de Fr\'echet induite 
sur leur cohomologie est pr\'eserv\'ee par ce quasi isomorphisme puisqu'il 
est continu (preuve laiss\'ee au lecteur) et coincide avec la structure
de Fr\'echet canonique sur $\Gamma(X^0,R^{\hidot}\tilde f _*{\cal 
G}^{\hidot})$.

Les structures de Fr\'echet de $\Gamma(X^0,R^{\hidot}\tilde f _*{\cal 
G}^{\hidot})$ et des espaces de cochaines de $A^{\hidot}$
(resp. de $B^{\hidot}$)  peuvent se laisser d\'ecrire par des familles 
d\'enombrables de
semi normes $\{ \|. \|_n\}_{n \in \mathbb N}$ de fa\c con que toutes 
les morphismes introduits ci dessus soient continus en norme $\|. \|_n$.

Soit $q\in {\mathbb Z}$. L'espace $A^{q}_2$ est par d\'efinition isomorphe 
\`a l'espace des suites
$(a^q_{\gamma})_{\gamma \in \Gamma/ \Sigma}$ avec $a^q_{\gamma}\in A^q$ 
telles que 
pour tout $n$, $\sum_{\gamma} \|a^q_{\gamma} \|_n^2 <\infty$. Comme il y a 
une description similaire pour $B^q_2$, l'inclusion $A^{\hidot}_2\subset 
B^{\hidot}_2$
r\'esulte de la continuit\'e de l'inclusion $A^{\hidot}\to B^{\hidot}$.

$A^{\hidot}$ \'etant strict,  appliquant le th\'eor\`eme de l'application 
ouverte pour les espaces de Fr\'echet \`a la surjection d'espaces de 
Fr\'echet 
$d:A^{q-1} \to dA^{q-1}$, on d\'eduit le complexe $A^{\hidot}_2$ est 
strict et que sa cohomologie se laisse d\'ecrire comme  
l'espace des suites
$(h^q_{\gamma})_{\gamma \in \Gamma/ \Sigma}$ avec $h^q_{\gamma}\in 
H^q(A^{\hidot})$ telles que 
pour tout $n$, $\sum_{\gamma} \|h^q_{\gamma} \|_n^2 <\infty$. 
La m\^eme description valant pour $B^{\hidot}_2$, il suit que 
l'inclusion
$A^{\hidot}_2\subset B^{\hidot}_2$
est un quasi isomorphisme. \qed

\paragraph{Suite spectrale de Leray}

\begin{prop}
Soit $\tilde X$ un $\Gamma$-espace complexe.
Soit $ {\cal F} \in C_{\bar\Gamma,\chi} (\tilde X)$. Soit $\tilde f:\tilde 
Y \to \tilde X$ un morphisme propre \'equivariant.

Il existe une suite spectrale de premier cadran
 $(E^{s,t}_r,d_r)_{r\ge 0}$ dans $Mod_{W^*(\bar \Gamma,chi)}$ aboutissant
\`a $H_{(2)}^*(\tilde Y,{\cal F})$ et  de terme $E_2^{r,s} 
=H_{(2)}^r(\tilde X,R^sf_*{\cal F})$.

Si  $\tilde X$ est cocompact, la suite spectrale $(E^{s,t}_r,d_r)_{r\ge 
0}$ peut \^ etre obtenue d\`es le terme $E_2^{r,s}$ \`a partir d'une suite 
spectrale dans 
$E_f(W^*(\bar\Gamma,\chi))$ par application du foncteur d'oubli 
\ref{isomalg}.
\end{prop}
{\bf Preuve:}
 La suite spectrale de Leray pour l'application continue $f:Y\to X$ et le 
faisceau $l^2\pi_*{\mathcal F}$ a pour terme $E_2$
 $$E^{r,s}_2=H^{r}(X, R^sf_* l^2\pi_*{\cal F})$$
 que \ref{tautologie} identifie canoniquement \`a 
$H^r( X, l^2\pi_* R^s\tilde f_*{\cal F})$. D'o\`u le premier point.
Le deuxi\`eme point est laiss\'e au lecteur (il s'agit d'une adaptation
facile de
la preuve du th\'eor\`eme \ref{tfini3}).\qed

\begin{coro}\label{devi}
Soit $\tilde X$ un $\Gamma$-espace complexe cocompact.
Soit $ {\cal F} \in C_{\bar\Gamma,\chi} (\tilde X)$. Soit $\tilde f:\tilde 
Y \to \tilde X$ un morphisme propre \'equivariant.
$$\sum_{n} (-1)^n \dim_{\bar\Gamma} H^n_{2} (\tilde Y,{\cal F})=
\sum_{r,s} (-1)^{r+s} \dim_{\bar\Gamma} H^r_{2} (\tilde Y,R^s\tilde 
f_*{\cal F})$$
\end{coro}

\subsection{Formules de Riemann-Roch par d\'esingularisation }

Le th\'eor\`eme de  r\'esolution des singularit\'es de Hironaka \cite{Hiro}
permet de ramener le calcul de la caract\'estique d'Euler 
d'un faisceau analytique coh\'erent sur un espace complexe au calcul de la 
caract\'eristique d'Euler de plusieurs faisceaux localement libres sur des 
vari\'et\'es complexes via la suite spectrale de Leray. Dans cet esprit, 
nous obtenons le:

\begin{theo}
Soit $\cal F$ un faisceau analytique coh\'erent sur l'espace 
complexe compact $X$. Soit $\pi:\tilde X \to X$ un rev\^etement galoisien
de groupe $\Gamma$. 
$$
\sum_{n} (-1)^n h^n_2(\tilde X,\pi^*{\cal F})=\sum_n(-1)^n h^n(X,{\cal F})
$$
\end{theo}
{\bf Preuve:} Par r\'ecurrence sur $\dim(X)$. En dimension $0$, ceci 
r\'esulte des propri\'et\'es ordinaires de la dimension de Von Neumann.  
Utilisant la filtration
$I$-adique ($I$ d\'esigne le radical de $X$), on se ram\'ene au cas o\`u $X$
est r\'eduit. On peut supposer que le support du faisceau ${\cal F}$ est 
$X$, sinon  l'utilisation de  \ref{devi} pour l'inclusion dans $\tilde X$
du $\Gamma$-espace complexe d\'efini par l'id\'eal annulateur de ${\cal F}$
permet de conclure.
Utilisant \ref{devi} pour une d\'esingularisation \'equivariante $\phi$,
on peut se ramener  au cas o\`u $X$ est lisse
puisque, pour $q\ge 1$ $R^q \phi_*\phi^*{\cal F}$ est support\'e en 
dimension $<\dim(X)$ et qu'
on a une suite exacte $0\to {\cal F}\to \phi_*\phi^*{\cal F} \to {\cal 
N}\to 0$, ${\cal N}$ \'etant support\'e en dimension $<\dim(X)$.
De m\^eme, utilisant la suite exacte
$0\to{\cal T}\to {\cal F} \to {\cal K}\to 0$ o\`u ${\mathcal T}$ est le 
sous-faisceau 
de torsion maximal, on se ram\'ene au cas o\`u ${\cal F}$ est sans torsion
sur la vari\'et\'e lisse $X$. On peut trouver $f:\hat X \to X$ un 
\'eclatement \'equivariant
 de $X$ tel que ${\cal V}=f^{-1} {\cal F} . O_{\hat X}$ est localement 
libre et $\hat X$ est lisse.
Le conoyau de l'injection canonique ${\cal F} \to f_*{\cal V}$
et les $R^qf_{*}{\cal V}$ sont
support\'es en dimension $<\dim(X)$, ce qui ram\'ene au cas d'un faisceau 
localement libre sur une vari\'et\'e, trait\'e par  \cite{Ati}. \qed

Le principe de d\'evissage \ref{devi} permet en principe de mener \`a bien 
le calcul
 de la caract\'eristique d'Euler $L^2$ d'un faisceau 
analytique coh\'erent projectivement p\'eriodique sur un $\Gamma$-espace 
complexe
propre cocompact en se ramenant au th\'eor\`eme d'indice $L_2$ d'Atiyah 
\cite{Ati}, \'eventuellement sous la forme plus g\'en\'erale 
utilis\'ee sous l'appelation \lq Vafa-Witten twisting trick\rq \ dans
\cite{Gro2} (voir \cite{Eys2}, pp. 189-196 pour une ex\'eg\'ese).

Une  suite spectrale de Leray en cohomologie $L_p$
existe tr\`es probablement et permet \'egalement de ramener l'\'etude de
la classe de $K$-th\'eorie  associ\'ee
\`a la somme altern\'ee des cohomologies $L_p$ d'un faisceau
analytique coh\'erent p\'eriodique
g\'en\'eral au cas d'un fibr\'e sur une base lisse. 
La plus petite cat\'egorie ab\'elienne dans laquelle vit la cohomologie 
$L_p$ avec $p\not=2$ \'etant loin d'\^etre aussi bien comprise que
$E_f(W^*_r(\Gamma))$, nous ne poursuivrons pas plus loin cette voie
puisqu'elle ne parait pas pouvoir amener
des applications substantielles.

\noindent
Philippe Eyssidieux, \\
Laboratoire Emile Picard, CNRS UMR 5580, UFR MIG \\
Universit\'e Paul Sabatier, 118 route de Narbonne \\
31062 Toulouse Cedex 4, France \\
e-mail: eyssi@picard.ups-tlse.fr


\begin{thebibliography}{99}
\bibitem{ABR} {D. Arapura, P. Bressler, M. Ramachandran}  {\em On the
fundamental group of a compact K\"ahler manifold}, Duke Math.
Journal {\bf 68} (1992) 477-488
\bibitem{Ati}{  M. Atiyah  }
{\em Elliptic operators,
 discrete groups and Von Neumann algebras}
 Soc. Math.Fr. Asterisque {\bf 32-33}, 1976, 43-72
 \bibitem{Cam}{F. Campana  }{\em Remarques sur le
rev\^etement universel
des vari\'et\'e k\"ahleriennes}, Bull. Soc. Math. Fr. {\bf
122} (1994), 255-284
\bibitem{Cam2} F. Campana {\em Fundamental group and
positivity of cotangent bundles of compact K\"ahler
manifolds}, J. Alg. Geo. {\bf 4} (1995), 487-502.
\bibitem{Cam3}{  F.  Campana  et J.P. Demailly }
Communication personnelle
\bibitem{ChGr} { J. Cheeger, M. Gromov}, {\em $L_2$-cohomology and
group cohomology}, Topology {\bf 25} (1986), 189-215
\bibitem{Dem}{ J.P. Demailly }{\em $L^2$ vanishing 
theorems for
Positive Line Bundles and Adjunction Theory}, in 
{\em Transcendental
Methods in Algebraic Geometry}, \'edit\'e par
F. Catanese et C. Ciliberto, LNM {\bf 1646} (1996), 
Springer
\bibitem{Dod}{ J. Dodziuk} {\em De Rham-Hodge theory for
$L^2$-cohomology of infinte coverings}, Topology {\bf 16} (1977),
157-165
\bibitem{Eys2} {P. Eyssidieux} {\em La caract\'eristique d'Euler du 
complexe de Gauss-Manin},
  Journal f\"ur die reine und
 angewandte Mathematik {\bf 490} (1997) 155-212
 \bibitem{Eys4} {P. Eyssidieux }{\em Th\'eorie de l'adjonction $L_2$
 sur le rev\^etement universel},
 Pr\'epublication 96 du Laboratoire
 Emile Picard, mai 1997.
 \bibitem{Eys5}{P. Eyssidieux }{\em Syst\`emes lin\'eaires adjoints 
$L_2$}, manuscrit 
soumis pour publication
\bibitem{Farb1} { M.S. Farber} {\em Homological
algebra of Novikov-Shubin invariants
and Morse inequalities}, GAFA, Vol. {\bf 6},
$n^o$ {\bf 4}, 1996
\bibitem{Farb2} M.S. Farber {\em Von Neumann Categories and
extended $L^2$ cohomology}, Pr\'epublication IHES/M/96/66, 1996
\bibitem{FoKn} O. Forster et K. Knorr {\em Ein Beweis des Grauertschen
Bildgarbensatzes nach Ideen von B. Malgrange}, Manuscripta Math. {\bf 5} 
(1971),
19-44 
\bibitem{Freyd} P. Freyd {\em  Abelian Categories}, Harper
and Row, 1964
\bibitem{GR} H. Grauert et R. Remmert {\em Coherent
Analytic
 Sheaves}, Grundlehren der mathemathischen Wissenschaften,
{\bf 265}, Springer, 1984
\bibitem{GR2} H. Grauert et R. Remmert {\em  Komplexe
R\"aume}, Math. Annalen {\bf 136}, 245-318 (1958).
\bibitem{GH} P. Griffiths, J. Harris {\em  Principles of Algebraic
Geometry}, John Wiley and Sons, 1978
\bibitem{Gro1} M. Gromov {\em  Sur le groupe fondamental d'une
vari\'et\'e k\"ahlerienne}, C.R. Acad. Sci. Paris {\bf 308} S\'erie 1
(1989), 67-70
\bibitem{Gro2}  
 M. Gromov, {\em K{\"a}hler hyperbolicity and
$L^2$-Hodge
theory}, J. Diff. Geom. {\bf 33}, 1991, 263-291
\bibitem{GHS} M.Gromov, G.Henkin, M.Shubin {\em  $L^2$
holomorphic functions on pseudoconvex coverings},
Pr\'epublication de l'IHES, IHES/M/95/58
\bibitem{GS} M.Gromov , M.Shubin {\em  Von Neumann spectra near
zero},
GAFA {\bf 1} (1991), 375-404
\bibitem{GS2}M.Gromov , M.Shubin {\em Near-cohomology of Hilbert
complexes and topology of non-simply connected manifolds},
Ast\'erisque {\bf 210} (1992), 283-294
\bibitem{Har} R. Hartshorne {\em Algebraic Geometry}, GTM {\bf 52} (1997), 
Springer
\bibitem{Ha2} R. Hartshorne {\em Residues and duality, Lecture Notes in 
Math} {\bf
20}, Springer, (1966)
\bibitem{Hiro} H. Hironaka {\em R\'esolution of 
singularities over a field of characteristic zero}, Ann. of 
Math.,{\bf 79} (1964), 109-326
\bibitem{JZ}J. Jost et K. Zuo {\em  Vanishing theorems for
$L^2$-cohomology on infinite coverings
of Compact K\"ahler manifolds and applications in algebraic
geometry}, Communication personnelle, Jan. 21, 1997.
\bibitem{Kol1}{Koll\'ar  J. } {\em Shafarevich maps and
plurigenera of algebraic varieties},
Inv. Math.{\bf 113} (1993), 177-215
\bibitem{Kol2}{ Koll\'ar J. }  {\em Shafarevich maps
and Automorphic Forms}, Princeton University Press (1995).
\bibitem{Lan} S. Lang {\em Algebra}, Addison-Wesley, $2^{nd}$ edition (1984)
\bibitem{LL} J. Lott et W. L\"uck {\em $L^2$-topological invariants
of $3$-manifolds}, Inv. Math {\bf 120}, 1995, 15-60
\bibitem{Luc1}W. L\"uck {\em 
 Approximating $L^2$-invariants by their finite dimensionnal counterparts}, 
GAFA
 {\bf 4}, 1994, 455-481
\bibitem{Luc2} W. L\"uck {\em Hilbert modules and modules over
finite Von Neumann algebras and applications to $L^2$ invariants}, 
Pr\'epublication (1996)
\bibitem{NR} T. Napier et M. Ramachandran {\em Structure theorems for 
complete K\"ahler 
manifolds and applications to Lefschetz type theorems}, GAFA {\bf 5} 
(1995), 809-851
\bibitem{NS} S. Novikov  et M.Shubin {\em Morse
inequalities and Von Neumann invariants of non simply
connected manifolds}, Sov. Math. Dok. {\bf 34} (1987), 79-82
\bibitem{NS1} S. Novikov  et M.Shubin, {\em Morse
inequalities  and Von Neumann $II_1$
factors}, Dokl. Akad. Nauk. {\bf 34} no. {\bf 1} (1986), 289-292
\bibitem{taka}   Takayama S. {\em Nonvanishing Theorems on
an algebraic Variety
with large fundamental group}, Pr\'epublication,
 Naruto University of Education,
(1997)
\bibitem{Verd} J.L. Verdier {\em Cat\'egories d\'eriv\'ees, \'etat 0}, 
LNM {\bf 569} (1977), Springer.
\bibitem{Yos} {K. Yosida} {\em Functionnal Analysis},
Grundlehren der mathematischen Wissenschaften, {\bf 123},
$6^{ieme}$ \'edition, Springer, 1980
\end{thebibliography}
\end{document}